 \documentclass{bull-l}

     \issueinfo{00}
     {}
      {}
     {2000}
   \copyrightinfo{2000}
     {American Mathematical Society}

               \usepackage{amssymb,amscd,verbatim,multicol}
\numberwithin{equation}{section}    
\makeatletter

\def\endoldequation{\eqno \@eqnnum 
$$\global\@ignoretrue}
\def\@eqnnum{\hbox to .01pt{}\rlap{\upshape
    \hskip -\displaywidth\tagform@\theequation}}
\makeatother

\let\cal\mathcal
\def\Ascr{{\cal A}}

\def\Cscr{{\cal C}}
\def\Dscr{{\cal D}}
\def\Escr{{\cal E}}
\def\Fscr{{\cal F}}
\def\Gscr{{\cal G}}
\def\Hscr{{\cal H}}

\def\Lscr{{\cal L}}
\def\Mscr{{\cal M}}
\def\Nscr{{\cal N}}
\def\Oscr{{\cal O}}
\def\OOscr{{\cal O}}
\def\Pscr{{\cal P}}
\def\Qscr{{\cal Q}}

\def\Uscr{{\cal U}}
\def\Vscr{{\cal V}}

\let\blb\mathbb
\def\CC{{\blb C}}

\def \PP{{\blb P}}

\def \ZZ{{\blb Z}}

\def \NN{{\blb N}}
\def \RR{{\blb R}}
\def \HH{{\blb H}}

\def\gr{\text{{\bf gr}-}}  
\def\hyphenn{\text{-}}
\def\Ga{\mathcal E}
\def\rr{\mathcal P}
\def\HH{\operatorname{H}}

\def\GKdim{\operatorname{GKdim}}
\def\GK{\GKdim}

\def\calL{\mathcal{L}}
\def\O{\mathcal O}
\def\RR{\mathcal R}
\def\pone{\mathbb P^1}
\def\too{\longrightarrow}
\def\Num{\operatorname{Num}}

\def\r{\mathbf r}

\def\Bimod{\operatorname{Bimod}}
\def\bimod{\operatorname{bimod}}
\def\Skl{\operatorname{Skl}}

\def\order{nice algebra}

\newcommand{\Proj}{\operatorname{Proj}}
\def\id{{\operatorname{id}}}
\def\Id{\operatorname{id}}
\def\Ob{\operatorname{Ob}}

\def\pr{\operatorname{pr}}

\def\Mod{\operatorname{Mod}}
\def\mod{\operatorname{mod}}
\def\Gr{\operatorname{Gr}}

\def\QGr{\operatorname{QGr}}
\def\qgr{\operatorname{qgr}}

\def\gr{\operatorname{gr}}

\def\Qch{\operatorname{Qch}}
\def\rep{\operatorname{rep}}
\def\wrep{\mathop{\widetilde{\rep}}}
\def\coh{\mathop{\text{\upshape{coh}}}}
\def\charact{\operatorname{char}}

\def\gr{\operatorname {gr}}
\def\Spec{\operatorname {Spec}}

\def\Ext{\operatorname {Ext}}
\def\Hom{\operatorname {Hom}}
\def\uHom{\operatorname {\mathcal{H}\mathit{om}}}
\def\End{\operatorname {End}}
\def\RHom{\operatorname {RHom}}

\def\Im{\operatorname {Im}}

\def\End{\operatorname {End}}

\def\Tails{\operatorname {Tails}}

\def\r{\rightarrow}

\DeclareMathOperator{\Tors}{Tors}
\DeclareMathOperator{\tors}{tors}

\DeclareMathOperator{\Aut}{Aut}

\def\dirlim{\mathop{\vtop{\baselineskip -100pt\lineskip
2pt\lineskiplimit 0pt
\setbox0\hbox{lim}\copy0\hbox to \wd0{\rightarrowfill}}}\limits}

\newtheorem*{convention}{Standard convention}

\newtheorem{lemma}{Lemma}[section]
\newtheorem{proposition}[lemma]{Proposition}
\newtheorem{theorem}[lemma]{Theorem}

\newtheorem{lemmas}{Lemma}[subsection]
\newtheorem{propositions}[lemmas]{Proposition}
\newtheorem{theorems}[lemmas]{Theorem}
\newtheorem{corollarys}[lemmas]{Corollary}

{

}

\newtheorem{thm}[lemma]{Theorem}

\newtheorem{prop}[lemma]{Proposition}
\newtheorem{cor}[lemma]{Corollary}

\renewcommand{\thequesns}

\theoremstyle{definition}

\newtheorem{definition}[lemma]{Definition}
\newtheorem{definitions}[lemmas]{Definition}
\newtheorem{examples}[lemmas]{Example}
\newtheorem{example}[lemma]{Example}

\newtheorem{remark}[lemma]{Remark}
\newtheorem{remarks}[lemmas]{Remark}
\newtheorem*{ack}{Acknowledgement}

\newdimen\uboxsep \uboxsep=1ex
\def\uboxn#1{\vtop to 0pt{\hrule height 0pt depth 0pt\vskip\uboxsep
\hbox to 0pt{\hss #1\hss}\vss}}

\def\uboxs#1{\vbox to 0pt{\vss\hbox to 0pt{\hss #1\hss}
\vskip\uboxsep\hrule height 0pt depth 0pt}}

\begin{document}
\title{Noncommutative curves and noncommutative surfaces}
\author{J. T. Stafford}
\address{Department of Mathematics, University of Michigan, Ann Arbor,
MI 48109, USA.} 
\email{jts@math.lsa.umich.edu}
\author{M. Van den Bergh}
 \address{Departement WNI,  Limburgs Universitair Centrum, 
 3590 Diepenbeek, Belgium.}
  \email{vdbergh@luc.ac.be}
  \thanks{The  first author
 was supported in part by an NSF
 grant}
\thanks{The second author is a senior researcher at the FWO}

\begin{abstract} 

In this survey article we   describe some geometric results in the theory of
noncommutative rings and, more generally, in the theory of abelian categories.

Roughly speaking and by analogy with the commutative situation, the
category of graded modules modulo torsion over a noncommutative graded
ring of quadratic, respectively cubic growth should be thought of as
the noncommutative analogue of a projective curve, respectively
surface. This intuition has lead to a remarkable number of nontrivial
insights and results in noncommutative algebra. Indeed, the problem of
classifying noncommutative curves (and noncommutative graded rings of
quadratic growth) can be regarded as settled. Despite the fact that no
classification of noncommutative surfaces is in sight, a rich body of
nontrivial examples and techniques, including blowing up and down, has
been developed. \end{abstract}

\maketitle
\tableofcontents

\section{Introduction}
\label{sec1}

It has always been clear that noncommutative rings are much more complicated
than commutative ones. In order to quantify this feeling Mike Artin and Bill
Schelter  started a project about 15  years ago which aimed at completely
classifying the objects in  at least some non-trivial classes of noncommutative
rings, the major class being the rings now called Artin-Schelter regular rings
of dimension three (see Section~\ref{secasregular}).  Algebraic geometry does
not generalize na\"{\i}vely  to noncommutative algebra;  yet, remarkably, one
can still carry across much of the intuition and technique. One can think of
these techniques as constituting a kind of noncommutative algebraic geometry.
The present article constitutes an introduction to this field.

This theory stays close in spirit to the classical theory of algebraic geometry
as described in textbooks such as \cite{H,Sh}. There are many other
approaches, some of which will be mentioned at the end of the introduction. Our
guiding theme will be to classify or at least understand rings of small
dimension, just as one natural approach to algebraic geometry would be to
emphasize  curves and surfaces. As in the classical case, this  will also
force us to use more abstract notions and these will be introduced as needed. 
Curiously, it seems that the most substantial applications of this principle
occur with projective rather than affine geometry. This is perhaps due to the
fact that, as localization is rarely available in the noncommutative universe,
most of the affine techniques do not work, whereas the more global and
categorical approaches to projective geometry can be generalized.
Thus we will  restrict our attention to graded rings.

To begin with think about the commutative case. We assume throughout that
$k$ is an algebraically closed field and that all rings are $k$-algebras.
Given a finitely generated, commutative graded $k$-algebra  
 $C=k[x_1,\ldots,x_{n}]/I$ that is generated in degree one,  
the zeroes of $I$ yield an associated projective
scheme $X=\Proj(C)\subseteq \PP^{n-1}$. Conversely, given $X$ one can recover
$C$ up to a finite dimensional vector space by defining $C_n$ to be the
functions on $X$ with a pole of order at most $n$ at infinity.  Crucially, the
modules over $X$ and $C$ are closely related; the category $\coh(X)$ of
coherent sheaves on $X$ is equivalent to the quotient category $\qgr(C)$
of finitely generated,  graded $C$-modules modulo those of finite length (see
Theorem~\ref{serre} for the precise statement).  In what follows we will
regard $\coh(X)$ rather than $X$ as the scheme, since this is what generalizes.
It can also  be  justified by
 Rosenberg's reconstruction theorem \cite{rosenberg1}.

Now consider a finitely generated, noncommutative graded $k$-algebra $R$.
Obviously $R$ can be written as a factor ring $T/I$ of a free
algebra $T$, 
 but this is not very useful  since such a
presentation gives virtually no insight into the structure of $R$ and
certainly does not seem to lead to any geometric insight.  

One therefore needs a more geometric way of studying such algebras.
Fortunately  the categorical aspects of the commutative theory do generalize
and one of the major themes of this article will be to work not just with $R$
but also with  
$\qgr(R)$, where these more geometric techniques apply. Crucially, one can go
back and forth:
\begin{itemize}
\item {\it {\em (}a{\em)} One can determine when
a category $\Cscr$ is isomorphic to $\qgr(R)$ for a noetherian ring $R$
{\em (}combine Theorems~\ref{AZ} and  \ref{Gabriel}{\em)}.

\noindent
 {\em (}b{\em)} For a large class of algebras 
 one can recover $R$, at least up to a finite dimensional vector
space, from data in $\qgr(R)$
{\em (}see Theorem~\ref{AZ2}{\em)}.

\noindent
 {\em (}c{\em)} Simple conditions on $R$ are known for
$\qgr(R)$ to satisfy geometric conditions like Serre duality {\em (}see
Section~\ref{secmorecomments}{\em)}.}  
\end{itemize}
 The most important special case of part (a) and 
the fundamental construction
 in the subject involves a noncommutative analogue of the homogeneous coordinate
 ring of a projective scheme $X$. This is  called the twisted 
homogeneous coordinate ring of $X$ and is defined as follows:
Given an automorphism $\sigma$ of 
 $X$ and an invertible sheaf $\Lscr$, the associated twisted homogeneous
coordinate ring is defined to be the ring 
$$B(X,\Lscr,\sigma)=\bigoplus_{i\geq
0} \HH^0(X, \Lscr\otimes \sigma^* \Lscr \otimes\cdots\otimes 
(\sigma^{n-1})^*\Lscr),$$ with a natural multiplication (see
Section~\ref{sec2.2}). Here, $\tau^*\Lscr$ denotes the pull-back of $\Lscr$
along the automorphism $\tau$. 
 Despite the fact that $B(X,\Lscr,\sigma)$  is  far from
commutative, one still has  $\qgr(B(X,\Lscr,\sigma))\simeq \coh(X)$ when
$\Lscr$ satisfies the appropriate form of ampleness. For example, the two most
basic noncommutative graded rings of quadratic growth are  the quantum (affine)
planes  $$U_\infty=k\{x,y\}/(xy-yx-x^2) \qquad\text{and}\qquad
U_q=k\{x,y\}/(xy-qyx), \quad\text{for } q\in k^*.$$   These can both be
written as $B(\PP^1, \Oscr_{\PP^1}(1), \sigma)$, where $\sigma$ is defined,
respectively,  by   $\sigma(u)=u+1$ and $\sigma(u)=qu$ on $k(\PP^1)=k(u)$. For
both rings $U$ one therefore has $\qgr(U)\cong\coh(\PP^1)$. 
 
Twisted homogeneous coordinate rings
seem to be pervasive
in the  geometrical approach to noncommutative rings being promoted in
this survey.
They  are  fundamental in the classification
of noncommutative curves. Perhaps more surprisingly, 
many graded rings $R$ come
equipped with a canonical map to a twisted homogeneous
 coordinate
ring $B$ (see Section~\ref{secasregular}) 
and 
that coordinate ring then forms the starting point in the 
study of the original ring $R$ (see
Sections~\ref{secasregular} and \ref{sec12}). 

Let us now  consider noncommutative curves. The  
first   question is to find  the appropriate definition of such an object.
  There are two possible answers, both of which play a
r\^ole in this subject. The first is algebraic and asserts that, just
as in the commutative case, {\it a curve is the abelian category
$\qgr(R)$ associated to a graded ring $R=\bigoplus_{i\geq 0} R_i$ with
quadratic growth; thus $\dim_k(R_i)$ grows linearly.}  The
rings $U_\infty$ and $U_q$ would be typical examples. The second
answer is geometric and asserts that {\it a smooth noncommutative curve
should be an abelian category 
$\Cscr$ with the basic properties of
$\coh(E) $ for a smooth projective curve $E$.} 
These analogies work rather well:
\begin{itemize}
\item{\it let $R$ be a graded domain of quadratic growth, and assume that
$R$ is generated by $R_1$ as
an algebra. Then, up to a finite dimensional
vector space, $R\cong 
B(E,\Lscr, \sigma)$, for some curve $E$. 
Moreover, the corresponding categories $\coh(E)$ and 
$\qgr(R)$ are  equivalent {\em(}see Theorem~\ref{thm3.3} and 
Corollary~\ref{cor3.4}{\em)}.}
\end{itemize}
This result allows one to give a detailed description of the structure of $R$;
for example, if $\sigma$ has finite order then $R$ will be  finitely generated 
as a  module over its centre, while if $\sigma$ has infinite order then 
$R$ has a faithful irreducible representation 
(see Corollary~\ref{cor3.5}).
When $R$ is not generated in degree one, the 
situation is somewhat more subtle although, under  mild 
assumptions, $R$ can still be classified in terms of 
a twisted coordinate ring $B(E,\Lscr,\sigma)$, and the categories 
$\coh(E)$ and 
$\qgr(R)$ will still be equivalent. See Sections~\ref{sec4} and \ref{sec5}
for the details.

  On the other hand, if one considers a category $\Cscr$ with the basic
  properties of $\coh(E)$ then:
  \begin{itemize}
  \item{\it such an abelian category
 $\Cscr$ is equivalent to $\qgr(R)$ for
the appropriate prime ring $R$ of quadratic growth
 {\em(}see Section~\ref{secy1}{\em)}.}
\end{itemize}
  These results give a precise formulation of the intuitive picture of
noncommutative curves, and we can regard this topic as being
essentially settled. 

The second half of the survey considers noncommutative surfaces.
Despite the fact that no form of classification is in sight, a rich body of
nontrivial examples and techniques has been  developed that should be useful,
not only for studying these surfaces but  possibly also for more general
algebraic questions.

For the purposes of this introduction, define a 
noncommutative irreducible surface to be $\qgr(R)$ for   a graded domain $R$
of cubic growth, as would be the case for commutative surfaces (see 
 Section~\ref{secnoncom-surface} for a more detailed discussion of 
 the appropriate definition).
In the commutative case, which is briefly reviewed in 
Section~\ref{seccomsur},
the simplest example is the projective plane 
$\PP^2 =\Proj (k[x,y,z])$.  As we
show (see Section~\ref{secncpq}) there are strong arguments for  saying that:
\begin{itemize}
\item{\it 
the noncommutative projective planes are $\qgr(A)$, for  the so-called
Artin-Schelter regular algebras $A$ of dimension $3$ with Hilbert series
$(1-t)^{-3}$ {\em(}see Section~\ref{secasregular}{\em)}. 
These rings have been completely classified.    }
\end{itemize}
Even in dimension three, these Artin-Schelter algebras 
are more diverse and have a more complicated structure than
might be expected. The obvious examples are the analogues of the quantum
affine planes mentioned earlier; the best known one  being 
the $k$-algebra with generators $\{x_i : 1\leq i\leq 3\}$ which
``$q$-commute''  in the sense that $x_ix_j=q_{ij}x_jx_i$, for 
$i<j$ and $q_{ij}\in k^*$. However the generic example is the  
Sklyanin algebra $\Skl_3$
   (see Example~\ref{skl} et seq) and this is much more subtle. In particular, 
   it  has a
factor ring isomorphic to $B(E,\Lscr,\sigma)$, where $\sigma$ is an
automorphism of a smooth  elliptic curve $E$. Geometrically, this can be
interpreted as saying that the noncommutative projective plane $\qgr(\Skl_3)$
has an embedded commutative elliptic curve $\coh(E)$. The algebra $\Skl_3$ 
 is actually  determined by the data  $(E,\Lscr,\sigma)$
 (the analogous statement holds for all  Artin-Schelter algebras of dimension
 three)
  and so, even here, the
noncommutative space is strongly controlled by  a commutative variety
via a twisted coordinate ring. 

The standard classification of smooth projective surfaces has three
parts:
first one classifies surfaces up to 
birational equivalence\label{birat-index} (this just
means
that one classifies the fields of transcendence degree two over $k$);
secondly
one has Zariski's theorem that asserts that one can pass from any smooth
projective surface to a birationally equivalent one by successive
blowing up
and down; thirdly one needs to know the minimal models within each such
equivalence class. The minimal models for rational surfaces are $\PP^2$,
 $\PP^1\times \PP^1$ 
 and $\PP^1$-bundles over $\PP^1$.  As mentioned in the last
paragraph, we know the noncommutative analogues for  $\PP^2$ and 
\begin{itemize}
\item{\it 
the noncommutative analogues
of $\PP^1\times \PP^1$ and  $\PP^1$-bundles over $\PP^1$ can also be
determined {\em(}see Sections~\ref{secncpq} and \ref{sec11}{\em)}; }
\end{itemize}
Moreover:
\begin{itemize}
\item{\it 
There even exists a
noncommutative version of blowing up and down {\em(}see
Section~\ref{sec12}{\em)}.}
\end{itemize}
  This is considerably more subtle than the commutative analogue since we need
to blow up at one-sided rather than at two-sided ideals and the arguments
therefore have to be categorical rather than ring-theoretic.
Nevertheless, the
end result is surprisingly similar to the commutative version. We
illustrate this  by showing that:
\begin{itemize}
\item{\it 
 the standard constructions of the Cremona
transformation, cubic surfaces and quadric surfaces via blowing up/down 
generalize completely to noncommutative surfaces.}
\end{itemize} 
This can be seen by comparing the
noncommutative versions in  (\ref{noncom-crem}--\ref{noncom-cub}) 
 with their commutative analogues, as outlined in
(\ref{seccub}--\ref{secquad}).

As was mentioned earlier, there is no complete classification for
noncommutative surfaces and, in particular, there is as yet no analogue
of the
classification of fields of transcendence degree two nor of Zariski's
theorem.
This will hopefully form a fruitful area of research for the reader.

There are many different approaches to classical geometry and---at least as
yet---these become more diverse when one considers their noncommutative
(or quantum or deformed) analogues. As we have argued above, the 
theory of projective curves and surfaces generalizes surprisingly effectively
(see  \cite{Kap1} for an application to the study of noncommutative
instantons). One also has  the seminal work by A. Connes, on what
is essentially  noncommutative differential geometry. Similarly,  much of the
theory of quantum groups can be thought of as a noncommutative geometry,
more precisely as   a noncommutative analogue of the  
theory of algebraic groups (see, for example, \cite{Manin3}). 
Just as differential geometry and the theory of algebraic groups  
tend to be rather different 
from classical algebraic geometry, so their noncommutative counterparts are
distinct.

Other  variants of noncommutative geometry can be found in
\cite{BaKo,Bondal3,KoRo,Leb,Manin,rosenberg, VOV} (this is by no means an
exhaustive list) and we will end the introduction  by briefly discussing a few
of these approaches to noncommutative algebraic geometry.

Equating
noncommutative schemes with abelian categories, as we have done, has
drawbacks, most notably there is no good general notion  of vector bundle. In
the recent preprint by Kontsevich and Rosenberg \cite{KoRo} this
 defect is elegantly repaired by
introducing a suitable notion of covering (see also \cite{VW}).  
Unfortunately the covering formalism creates other
difficulties and it is not clear whether a noncommutative graded ring
should give rise to a noncommutative scheme in their sense.

 Many approaches to noncommutative algebraic geometry
 may eventually be incorporated
into a theory of deformations of orders over varieties. 
However, this first requires a thorough understanding of the geometry 
of such orders or, more generally, of rings satisfying a polynomial
identity. For  the latter the reader is
referred to \cite{VOV}. Of course, the fundamental property of 
rings satisfying a polynomial
identity is that their irreducible representations are finite dimensional,
and for a general theory of the moduli spaces of finite dimensional
representations of more general rings the reader is referred 
to \cite{Leb,Leb1}. Closely related is Kapranov's notion of a 
noncommutative ``thickening'' of a smooth commutative variety \cite{Kap}.
 Many of the questions studied in these papers are related to
Kontsevich's suggestion that many  moduli problems could be
reformulated as a study of the finite dimensional representations of 
infinite dimensional algebras.

Arguably \cite{Bondal3} one should generalize  further and work
directly with triangulated categories. However, due to a certain
flabbiness in their definition this can create serious technical
problems that are quite unrelated to the underlying geometric
intuition.  It is now clear that most of these problems can be solved
by considering DG-categories~\cite{Keller1}
 or more generally A${}_\infty$ categories~\cite{Keller}
  so it seems worthwhile to pursue this viewpoint
further.

\begin{ack}
We would like to thank Dennis Keeler and Haynes Miller
for many suggestions that improved the exposition of this article. 
\end{ack}
\subsection{Notation} \label{notation}
Here, we collect some  standard notation that will be in force
 throughout the paper.

 Fix once and for all a base field $k$, which will be algebraically closed
 unless we say otherwise.  
All objects and morphisms will be assumed to be $k$-linear. Thus all
rings will be $k$-algebras, all schemes will be $k$-schemes, all additive 
categories will be $k$-linear, etc.  If $A$ is a ring then $\Mod(A)$ will
denote the category of right $A$-modules. Similarly, if
$R=\bigoplus_{i\in \ZZ}R_i$ is a $\ZZ$-graded ring then 
$\Gr(R)$\label{Gr(R)-index} will
be the category of $\ZZ$-graded right $R$-modules. In this category
homomorphisms are of degree zero; thus if $M=\bigoplus_{i\in \ZZ} M_i$
and $N=\bigoplus_{i\in \ZZ} N_i \in \Gr(R)$, and $\theta\in
\Hom_{\Gr(R)}(M,N)$, then $\theta(M_i)\subseteq N_i$ for all $i$. We let
 $\Tors(R)$\label{Tors(R)-index} denote the full subcategory of
$\Gr(R)$ generated by the right-bounded modules; that is, by modules
$M=\bigoplus_{i\in \ZZ} M_i$ for which $M_i=0$ for $i\gg 0$.  Set
$\QGr(R)$\label{Qgr(R)-index} for the quotient category
$\Gr(R)/\Tors(R)$ (the formal definition of this can
be found, for example in \cite{Gabriel} but, roughly speaking, $\QGr(R)$ has
the same objects as  $\Gr(R)$ but objects in  $\Tors(R)$ become isomorphic
to $0$).   Given a graded module $M=\bigoplus M_i$ over a 
graded ring $R$,  define
the {\it shift} $M(r)$\label{shift-index} of $M$ to be the graded module defined by
$M(r)_i=M_{i+r}$.

For the most part we will work with
$\NN$-graded algebras $R=\bigoplus_{i\geq 0} R_i$ that are finitely
generated, with $R_0$ finite dimensional as a $k$-vector space. We
call these algebras {\it finitely graded}\label{fingraded-index} for short. 
A finitely graded ring $R$ will be called {\it connected
graded}\label{conngraded-index} if
$R_0=k$. The difference between these two concepts is largely technical;
in particular, if $R$ is a finitely graded domain then $R$ is connected
graded.  For an  $\NN$-graded algebra $R$
we write $R_+$ for the
two-sided ideal $\sum_{n>0} R_n$.

 We should emphasize that, although we are interested in
noncommutative analogues of geometric concepts, the unadorned terms ``scheme,'' 
``variety,'' etc, will always mean the classical (commutative)
objects.   
Let $X$ be a scheme. Then  $\Qch(X)$\label{Qch(X)-index}
 will denote the category of
quasi-coherent $\Oscr_X$-modules. If $X$ is noetherian then we
denote by $\coh(X)$\label{coh(X)-index}
 the category of coherent $\Oscr_X$-modules. In all
cases of practical interest   $\Qch(X)$ is a so-called
\emph{Grothendieck
  category}\label{grot-index} \cite{thomasson}.
  This is an abelian category with a generator and exact
filtered direct limits. Such a category automatically 
has products and
enough injectives \cite[Th\'eor\`eme 1.10.1]{Grot}. 
Since abelian categories can be embedded in module 
categories, many standard terms from module theory 
extend naturally to abelian
categories and, on the whole, their definitions will 
not be given explicitly.

 {\def\thelemma{(*)}
\begin{convention} \label{convent} Throughout the paper, if
$\mathrm{XYuvw}(\cdots)$  denotes an 
  abelian category then  we will denote by 
$\mathrm{xyuvw}(\cdots)$ the
full subcategory consisting of noetherian 
objects (an object is noetherian if it satisfies the ascending chain
  condition
on subobjects).
\end{convention} 
}
This convention fixes the
meaning of $\mod(A)$\label{gr(R)}, $\gr(R)$, $\tors(R)$ and $\qgr(R)$. Note that,
when $R$ is finitely graded,    
$\tors(R)$ consists of the finite dimensional graded $R$-modules.
As is true of $\Mod(R)$ versus  $\mod(R)$, it is 
sometimes more convenient
(and sometimes less)  to work  with 
 $\QGr(R)$ rather than $\qgr(R)$ or with  
$\Qch(X)$ rather than $\coh(X)$. The following theorem 
says that each member of these pairs determines the other, and so 
we are free to work in whichever is most convenient. 

This also works more abstractly. 
Recall that an abelian category is {\it noetherian}\label{noethcat-index}
if every  object is noetherian  and  the isomorphism classes
of objects form a set. 
Given such a noetherian abelian 
category $\Cscr$ we may embed $\Cscr$ in its closure $\widetilde{\Cscr}$
  under direct limits \cite{Gabriel}. Then $\widetilde{\Cscr}$ is 
a  Grothendieck category, that  is
{\it locally noetherian}\label{locnoeth-index}
 in the sense that it has a family
of noetherian generators. 
Since $\tilde{\Cscr}$  has enough injectives one can, for example,  
 define the derived functors
 of Hom groups from $\Cscr$ by working in~$\widetilde{\Cscr}$.

\begin{theorems}\label{Gabriel}\cite[Theorem II.4.1]{Gabriel}
 Let $\Cscr$ be a noetherian  abelian category and $\widetilde{\Cscr}$
its closure, as above. Then $\Cscr$ is the set of noetherian objects in 
 $\widetilde{\Cscr}$ and  these two categories determine each other up to
 natural equivalence.

If $\Cscr=\coh(X)$ for a noetherian scheme $X$, then
$\widetilde{\Cscr}=\Qch(X)$ and if $\Cscr=\qgr(R)$ for a right noetherian
 graded  ring $R$, then $\widetilde{\Cscr}=\QGr(R)$
\end{theorems}

Finally, the following ring-theoretic definitions 
will be used without comment in
the survey. 
A ring $R$ is {\it prime},\label{prime-index}
respectively {\it semiprime}, if the product of nonzero ideals is
nonzero, respectively the square of every nonzero ideal is nonzero.
An ideal $I$ of $R$ is {\it (semi)prime}\label{semiprime-index}
 if $R/I$ is (semi)prime.
The {\it prime radical}\label{radical-index} $N(R)$ of  $R$ is
defined to be the intersection of the prime ideals of $R$. 
A useful fact is that $R$ is semiprime if and only if $N(R)=0$. 
Given $R$-modules $M\subseteq N$, then $N$ is an 
{\it essential}\label{essent-index} extension of
$M$ if every nonzero submodule of $N$ intersects $M$ nontrivially.

\section{Geometric Constructions}
\label{sec2}
The purpose of this section is to clarify the connection between rings
and abelian categories in certain standard cases. The main result (see
Theorems~\ref{AZ} and \ref{AZ2}) determines which categories are equivalent to
the category $\QGr(R)$ over a right noetherian ring $R$.
This will lead us to the
notion of ``twisted homogeneous coordinate rings'' which will turn out
to be a basic tool  in the construction and 
classification of noncommutative curves. 
A very simple example 
 that may help to
illustrate these ideas is given in \eqref{simple-eg}.

\subsection{Affine and projective schemes}
\label{sec2.1}
Although our main interest will be projective schemes, let us fix the
ideas 
by first
considering the affine case.
  Let $X$ be an affine scheme and put
$A=\Gamma(X,\Oscr_X)$. Then it is classical that $\Qch(X)$ is
equivalent to $\Mod(A)$. This  immediately begs the question which
Grothendieck categories can be written as $\Mod(A)$ for
 a possibly noncommutative ring $A$. The answer  is again classical:
\begin{propositions}\label{affineequiv} \cite[Example
X.4.2]{stenstrom}
Let $\Cscr$ be a Grothendieck category with a  projective
generator $\OOscr$
 and assume that $\OOscr$ is small; that is, assume that 
 $\Hom(\OOscr,-)$ commutes with
direct sums. Then $\Cscr\simeq\Mod(A)$, for $A=\End(\OOscr)$.
\end{propositions}

Now  consider the projective case. Let $X$ be a  projective scheme
with a line bundle $\Lscr$. Then the  \emph{homogeneous
  coordinate} ring $B=B(X,\Lscr)$ 
  associated to $(X,\Lscr)$ is defined by the formula
$
B=\bigoplus_{n\in \NN} \Gamma(X,\Lscr^n)
$
with obvious multiplication.  Similarly if $\Mscr$ is a quasi-coherent
sheaf on $X$, then 
$
\Gamma_h(\Mscr)=\bigoplus_{n\in \NN} \Gamma(X,\Mscr\otimes\Lscr^n)
$\label{gamma-h-index}
defines a  graded $B$ module. Thus 
  the composition of $\Gamma_h$ with the natural projection from
$\Gr(B)$ to
  $\QGr(B)$ yields a functor $\overline{\Gamma}_h : 
  \Qch(X)\r \QGr(B)$. This functor only works well when $\Lscr$ is 
  {\it ample}\label{amplediv-index} in the sense that, if $\Fscr\in \coh(X)$,
 then $\Fscr\otimes
\Lscr^{\otimes n}$ is generated by global sections
for all $n\gg 0$ 
\cite[p.153]{H}.  In this case one has the following fundamental
result of Serre:

\begin{theorems}\label{serre} \cite[Proposition 7.8]{Se}
{\em(i)} Let $\Lscr$ be an ample line bundle on
 a projective scheme $X$. Then the functor $\overline{ \Gamma}_h(-)$
defines an equivalence of categories between $\Qch(X)$ and $\QGr(B)$.

\noindent {\em (ii)} Conversely, if $R$ is a commutative 
 connected graded
$k$-algebra that is generated by $R_1$ as a $k$-algebra, 
 then there
exists a line bundle $\Lscr$ over $X=\Proj(R)$ such that $R=B(X,\Lscr)$, up to
a finite dimensional vector space. Once again, $\QGr(R)\simeq \Qch(X)$.
\end{theorems}
It is therefore 
 natural to  think of Grothendieck categories of the form
$\QGr(B)$ with $B$ a (possibly noncommutative) graded ring as
\emph{noncommutative projective schemes}.
      Of course, as always, this sort of
intuition is only worthwhile if it leads to something useful, but we
hope that
this survey shows that it does!

\begin{remarks}\label{serre2}
If $R$ is a connected graded  commutative ring that 
is not generated in degree one,
then, typically, $\QGr(R)\not\simeq\Qch(\Proj(R))$ and so 
   Theorem~\ref{serre}(ii) fails. For example, if $R$ is a
weighted
  polynomial ring, in the sense that the variables are given degrees
other than
  one, 
   then $\QGr(R)$ clearly has finite homological
    dimension whereas $\Proj(R)$ is usually singular.
    \end{remarks}
    
The first obvious question is to determine which  Grothendieck
categories $\Cscr$ are of the form  $\QGr(B)$.  
As we next explain, this has been solved 
by Artin and Zhang \cite{AZ}. 
See \cite{V} for related results of Verevkin. 

In Serre's Theorem, one should regard the equivalence as being
controlled by
three pieces of geometric data: the category $\Cscr=\Qch(X)$,  the
object
$\Oscr_{X}$ and the autoequivalence of $\Cscr$ given by tensoring with
$\Lscr$. 
 This has a natural generalization. 
Let $(\Cscr,\OOscr,s)$ be a triple consisting of a Grothendieck category
$\Cscr$, an object $\OOscr$ in $\Cscr$ and an
autoequivalence $s$ on $\Cscr$.
Associated to
this data one can define a homogeneous coordinate ring
\begin{equation}
\label{factory}
\Gamma_h(\Cscr,\OOscr,s)=\bigoplus_{n\ge 0} \Hom(\OOscr,s^n \OOscr) 
\end{equation}
with multiplication $a.b=s^n(a)b$ for $a\in \Hom(\OOscr,s^m \OOscr)$,
$b\in
\Hom(\OOscr,s^n\OOscr)$. 
To simplify the notation we write
$\OOscr(n)$ for $s^n\OOscr$.
\begin{definitions}
\label{ampleness}
Assume that $\Cscr$ is locally noetherian. 
Let
   $\OOscr\in\Ob(\Cscr)$ be a noetherian object and let $s$ be an
   autoequivalence of $\Cscr$.  Then    the pair $(\OOscr,s)$ is 
   called {\em ample} if the following
  conditions hold:
\begin{enumerate}
\item[(i)]
For every noetherian object $A\in\Cscr$ there are positive
integers $l_1,\ldots,l_p$ and an epimorphism from
$\bigoplus_{i=1}^p \OOscr(-l_i)$ to $A$.
\item[(ii)]
For every epimorphism   between noetherian objects $A\r B$  the induced
map $\Hom(\OOscr(-n),A)\r \Hom(\OOscr(-n),B)$ is surjective for $n\gg
0$.
\end{enumerate}
\end{definitions}
When $(\Cscr,\OOscr,s)=(\Qch(X),\Oscr_X,-\otimes\Lscr)$, this definition
of
ampleness does reduce to the standard one. Indeed, part (i) of the definition 
reduces to 
 the standard definition of an ample  sheaf, while part (ii) corresponds to 
  part
of the homological definition.  Moreover Serre's theorem
also generalizes, as is shown by the following special case of a result
of Artin and Zhang.
\begin{theorems} \label{AZ} \cite[Theorem 4.5]{AZ}
 Assume that $\dim_k\Hom(M,N)<\infty$ for all 
noetherian objects $M,N\in\Cscr$ and let $(\OOscr,s)$ be an ample pair.
Then  $ \Cscr\simeq \QGr(B) $ for $B=\Gamma_h(\Cscr,\OOscr,s)$. In addition
$B$ is right  noetherian.
\end{theorems}

In this theorem, the functor $\overline{\Gamma}_h
 : \Cscr \to\QGr(B)$ is provided by
an analogue  of the functor used in Theorem~\ref{serre}: 
For $\Mscr\in\Ob(\Cscr)$,
define $\overline{\Gamma}_h(\Mscr)$ to be the image in $\QGr(B)$ of 
$\Gamma_h(\Mscr)=\bigoplus \Hom_{n\in \NN}(\OOscr, s^n\Mscr)\in
\Gr(B)$.

There is also a converse to this theorem, but this requires an extra
hypothesis since  $B=\Gamma_h(\Cscr,\OOscr,s)$ 
satisfies the extra condition:
 \begin{equation}\label{chi}
\chi_1:\quad\dim_k \Ext_{B}^1(B/B_+, M)<\infty,
\end{equation} for all finitely
generated  graded $B$-modules $M$. If we allow
this condition, then the natural converse to Theorem~\ref{AZ} goes
through.  We let $\pi M$\label{pi-index}
 denote the image in $\QGr(R)$ of $M\in \gr(R)$.

\begin{theorems} \label{AZ2} \cite[Theorem 4.5]{AZ}
  Let $B$ be a finitely graded,  right noetherian $k$-algebra that
  satisfies $\chi_1$.  Set $(\Cscr,\OOscr,s)=(\QGr(B),\pi B, (+1))$,
where  $(+1)$ is the shift functor defined in \S\ref{notation}.
  Then  $(\pi B,(+1))$ is ample and
  $\dim_k\Hom_\Cscr(M,N)<\infty$ for noetherian objects
$M,N\in\Ob(\Cscr)$.
  Moreover, $B= \Gamma_h(\Cscr,\OOscr,s)$, up to a finite
  dimensional vector space.
\end{theorems}

For a more detailed analysis of this categorical approach
  we refer the
reader to \cite{AZ,AZ2}.

The intrinsic reason for requiring $\chi_1$ in this theorem 
is as follows: 
Let $M$ be a
finitely generated torsion-free $\ZZ$-graded  module 
 over a finitely graded
noetherian ring $B$. Then  \cite{AZ} shows that 
$\widetilde{M}=\Gamma_h(\pi M)$  is an $\NN$-graded 
essential  extension of $M$ by a torsion module. Thus, 
if $\widetilde{M}/M$ is
finite dimensional then $\dim_k\Ext^1_B(B/B_+, M)<\infty$. 
A significant step 
in the proof of Theorem~\ref{AZ2} is to prove the converse.

Unfortunately, not all noetherian connected graded rings satisfy
$\chi_1$.
The basic  counterexample \cite{Staf4} is provided by 
 \begin{equation}\label{non-chi}
 S=k + Uy\quad\subset \quad
U=k\{ x,y\}/(xy-yx-x^2).
\end{equation} 
 This ring is fundamental to many of the results in this paper
and so we should describe its properties in some detail. In the notation of
the introduction,  $U=U_{\infty}$  is  one of the two  quantum affine 
planes and one can construct a similar example using
 $U_q$. 
 As $U$ has basis $\{x^iy^j\}$, it follows 
that $S$ has basis $\{1, x^iy^j : i\geq 0, j\geq 1\}$. 
{}From the identity 
 $(x^{n-1}y)(xy) -(x^ny)y = -x^{n+1}y$,  
 it is easy to see that $S$  is generated by $y$ and $xy$. 
 However, notice that, although $U/S$ is infinite dimensional,
  $U\cdot S_+ = Uy\subset S$. Thus $U/S\subseteq \Ext^1_S(k, S)$ 
  and so $\chi_1$
  certainly fails.   Finally:
  \begin{equation}\label{s-is-noether}
  \text{$S$ is (right or left) noetherian if and only if char$\,k=0$.}
  \end{equation}
  The intrinsic reason why this holds will appear in  \S\ref{sec4},
  but it is easy to prove directly \cite[Theorem 2.3]{Staf4}. 
   Another amusing property of this ring is that the 
  presentation and Hilbert series of $S$ are independent of the
characteristic
  of $k$. Thus 
 \eqref{s-is-noether} implies that the noetherian
property fails
  to be preserved under reduction modulo $p$ for any prime $p$!

Theorem~\ref{AZ2} does not fail too badly for $S$, since the ring 
$\Gamma_h(\QGr(S),\pi S,(+1))$ obtained from that theorem is simply 
$U$---this should come as no surprise since the last paragraph shows that
$U\subseteq \widetilde{S}_S$ and it is easy to see that 
$\Gamma_h(\QGr(U),\pi U, (+1)) = U$.  
 However,  with more complicated examples one can
exhibit nastier behaviour.
 For example,  the polynomial ring $T=S[z]$ satisfies
  $\Gamma_h(\QGr(T),\pi T,(+1)) = T$, but Theorem~\ref{AZ2} must still
  fail, since it fails for the factor ring $S$.  What now goes wrong is
  that the shift functor $s=(+1)$ in $\QGr(T)$ is not ample; specifically,
 Definition~\ref{ampleness}(ii) fails for the map
  $T\twoheadrightarrow S$ (see \cite[~Theorem~2.10]{Staf4}).

The condition $\chi_1$, and the analogues $\chi_n$,
 obtained by replacing $\Ext^1$ by
$\Ext^n$ in its definition, are not well understood, although they do
hold for
most of the natural classes of algebras. They also seem to be
fundamental
to a deeper understanding of the homological properties of
noncommutative
graded rings.  See \S\ref{secmorecomments} for   further
results,
 questions  and references in this direction.


\section{Twisted homogeneous coordinate rings}
\label{sec2.2}
One may view \eqref{factory} as a factory for producing graded rings
and the interesting point about Theorem  \ref{AZ} is that fixing $\Cscr$ 
but varying
$(\OOscr,s)$ will produce many different graded rings $B$ with the same
$\QGr(B)$. As far as this survey is concerned, the most interesting examples
occur when $\Cscr = \Qch(X)$, for a projective scheme $X$.
We consider this case in detail in this section.
There are now two basic ways of obtaining an autoequivalence of
$\Cscr=\Qch(X)$; one can either tensor with a line bundle $\Lscr$
(as happens in the commutative case)
or apply $\sigma_\ast$ for $\sigma$ an automorphism of $X$. 
In the sequel we will simply combine them by
 considering the composition $s=\sigma_\ast(-\otimes
\Lscr)$. We will say that $\Lscr$ is \emph{$\sigma$-ample}\label{s-ample-index}
 if $(\Oscr_X,s)$ is
ample in the  sense of Definition \ref{ampleness}.

Following tradition we will write $B(X,\Lscr,\sigma)$ for
$\Gamma_h(\Qch(X), \Oscr_X,s)$ 
and call it the \emph{twisted homogeneous
coordinate ring} associated to the pair $(\Lscr,\sigma)$. 
If $\Mscr\in \Qch(X)$, then a trite
computation yields
\begin{equation}
\label{trite}
s^n(\Mscr)=\sigma_{\ast}^n(\Mscr\otimes \Lscr\otimes
\Lscr^\sigma\otimes\cdots \otimes \Lscr^{\sigma^{n-1}})
\end{equation}
where $\Lscr^{\tau}$ is a shorthand for the pullback
$\tau^\ast\Lscr$ along an automorphism $\tau$. 
 From \eqref{trite}
 one obtains $B(X,\Lscr,\sigma)=\bigoplus_n B_n$ where
$$B_n=\Gamma(X,\Lscr_n)\qquad {\rm with} \quad\Lscr_n=\Lscr\otimes
\Lscr^\sigma\otimes\cdots \otimes \Lscr^{\sigma^{n-1}}.$$  
The multiplication on $B=B(Y,\Lscr,\sigma)$ is defined in the obvious way:
We have a natural map
 \begin{align}\label{equ2.2}
B_n\otimes_k B_m  &\cong 
\HH^0({{X}}, \calL_n)\otimes_k \HH^0({{X}}, \calL_m)
\\ \nonumber
&\cong
\HH^0({{X}}, \calL_n)\otimes_k \HH^0({{X}}, \calL_m^{\sigma^n})
{\buildrel\phi\over\too}\HH^0({{X}}, \calL_{n+m}) = B_{n+m},
\end{align}
where the map $\phi$ is obtained by taking global sections of the natural 
homomorphism
$\calL_n\otimes\calL_m^{\sigma^n} \;{\buildrel{\sim}\over {\too}}\;\;
\calL_{n+m}.$

One reason why  twisted homogeneous coordinate rings are so important
 is that they frequently appear as factor rings of noncommutative graded
rings; indeed, they satisfy a universal property analogous to passing from a
noncommutative ring $A$ to its abelianization $A/[A,A]$. Since this
is technical to write out formally, we will refer the reader to
\S\ref{secasregular} where this process is illustrated by example, 
and to \cite{AZ2} for a more thorough treatment.

Of course, this discussion 
begs the question of precisely which line bundles are
$\sigma$-ample. A simple application of the
Riemann-Roch theorem shows that, 
\begin{equation}\label{rroch}
\text{\it if $X$ is a curve, then
 any ample invertible sheaf is $\sigma$-ample,}
 \end{equation}
and the converse holds for irreducible curves. However, for higher
dimensional 
varieties, the situation is more subtle and is described by
the following result, for which we need some notation.
Let $X$ be a projective scheme
and write  $A^1_{\Num}(X)$ for  the set of Cartier divisors of $X$
modulo 
numerical equivalence. 
Let $\sigma$ be an automorphism of $X$ and let $P_\sigma$ denote its
induced 
action on  $A^1_{\Num}(X)$. Since $A^1_{\Num}(X)$ is a finitely
generated  free
abelian group \cite[Remark 3, p.305]{Kl},
 we may represent $P_\sigma$ by a matrix, and we
call 
$P_\sigma$ {\it quasi-unipotent}\label{quasiunipot-index}
 if all the eigenvalues of this matrix
are
roots of unity.

\begin{theorem}\label{sigma-ample} \cite{Ke}
Let $X$ be a  projective scheme
and let $\sigma$ be an automorphism of $X$.
Then  $X$ has a $\sigma$-ample line bundle if and only if 
$P_\sigma$ is quasi-unipotent. If  $P_\sigma$ is quasi-unipotent, then 
all ample line bundles are $\sigma$-ample.
\end{theorem} 

One consequence of the proof of this result is that a line bundle
$\Lscr$ is
$\sigma$-ample if and only if it is $\sigma^{-1}$-ample.
This is significant since the opposite ring
$B(X,\Lscr,\sigma)^\circ\cong
B(X,\Lscr,\sigma^{-1})$.  Combined with 
Theorem~\ref{AZ}, this proves:

\begin{theorem}\label{AZZ} 
 Let $X$ be a projective scheme with an automorphism $\sigma$
and suppose that  $\Lscr$ is a $\sigma$-ample line bundle. 
Then  $B=B(X,\Lscr,\sigma)$ is
a left and right noetherian ring for which 
$\qgr(B) \simeq \coh(X)$.
\end{theorem}

 This result can be interpreted as saying that, although twisted
homogeneous coordinate rings are generally very noncommutative, they
behave from a geometric standpoint as if they are commutative. 

Checking the condition that $P_\sigma$ is quasi-unipotent is rather
easy in practice. For example,
 assume that $X$ is a   projective variety
with a canonical divisor that is either
  ample or anti-ample. Then $P_\sigma$ is quasi-unipotent for every
 automorphism  $\sigma$ \cite[Proposition 5.6]{Ke}. 
 It follows that for varieties like Grassmannians 
 all automorphisms must be quasi-unipotent. 
Furthermore it looks as if the condition for  $P_\sigma$ to be
quasi-unipotent might actually be topological. This is illustrated by
the following  result  due to  Mike
Artin. 

\begin{proposition}
 Assume that $X$ is a smooth projective surface over $\CC$. Then
  $P_\sigma$ is quasi-unipotent if and only if the action of $\sigma$
  is quasi-unipotent on $H^2(X,\CC)$.
\end{proposition}


\begin{example}\label{simple-eg}
Let us compute the  twisted homogeneous coordinate rings  in the
simplest nontrivial  case,  when
${{X}}$
is the projective line $\PP^1$ and $\calL=\O_{\pone}(1)$. 
Pick  a parameter $u$ for $\pone$, 
so that the standard affine cover $U=\pone \smallsetminus\infty$,
$V=\pone \smallsetminus 0$ has rings of regular functions
$\O(U)=k[u]$ and $\O(V)=k[u^{-1}]$. Identify $\O_{\pone}(1)$ with the
sheaf of
functions on $\pone$ that have at most a  simple 
pole at infinity; in other words, 
it is the  sub-sheaf
of $k(u)=k(\pone)$ generated by $\{1,u\}$.
Then a simple calculation shows that $\text{H}^0({{X}},\calL^{\otimes
n})$
is spanned by $\{1,u,\dots, u^n\}$ and that $B({{X}},\calL)\cong
k[x,y]$, where
$x=1 $ and $y=u$, thought of as elements of $B_1 =
\text{H}^0({{X}},\calL)$.

Now repeat this computation for $B(\pone, \O_{\pone}(1), \sigma)$  where
$\sigma$ is the automorphism defined by $\sigma(u)=u+1$. Here, we use
the same
symbol to denote both  an automorphism of the field of rational
functions 
$K=k({{X}})$ of a smooth projective curve ${{X}}$ and for the
automorphism this induces on ${{X}}$. 
Our convention is that  $\sigma$ acts on the right on $K$
and on
the left on ${{X}}$; thus $f^\sigma(y)=f(\sigma(y))$, for $f\in K$ and
$y\in
{{X}}$. Then  $\calL=\O_{\pone}(1)$ satisfies  $\calL^\sigma \cong
\calL$;
indeed under our concrete presentation of $\calL$,  one even has
$\calL^\sigma
= \calL$ and hence $\calL_n=\calL^{\otimes n}=\O(n)$. Therefore,
$B(\pone,\calL,\sigma)=B(\pone,\calL,\Id)$, as graded vector spaces.
However,
the multiplication is twisted: If we again set $x=1$ and  $y=u$, as
elements of
$B_1$, then the multiplication map \eqref{equ2.2} sends  $$y\otimes x
=u\otimes
1 \mapsto u\otimes 1^\sigma =u\otimes 1 \mapsto u \in \text{H}^0(\pone,
\calL_2).$$ On the other hand,  
$$x\otimes y =1\otimes u \mapsto 1\otimes
u^\sigma =1\otimes (u+1) \mapsto u+1 \in \text{H}^0(\pone, \calL_2).$$ 
Since
$x\otimes x\mapsto 1\in \text{H}^0(\pone, \calL_2)$,  we see that we
have the
relation $xy=yx+x^2$, in $B(\pone,\calL,\sigma)$.  As one might expect,
a
little more work shows that  
\begin{equation}\label{babyeg}
B(\pone,\calL,\sigma)\cong k\{ x,y\}/(xy-yx-x^2).
\end{equation}
 We leave it as an exercise to the reader to show that, if one uses the
 automorphism $\sigma:u\mapsto \lambda u$, for some $\lambda\in k^*$,
 then $B(\pone,\O_{\pone}(1),\sigma)$ is isomorphic to the other 
 quantum plane $k\{ x,y\}/(xy-\lambda yx).$
 \end{example}

\section{Algebras with Linear Growth.}\label{sec2a}

As was pointed out in the introduction, it is necessary to
have some type of dimension function on noncommutative rings. 
There is a multitude of
such functions in the noncommutative case, each with its own
advantages and disadvantages. The interested reader may look in
\cite{KL,MR,YZ1,Zhang2} for some relevant references.

As an example of a reasonably well-behaved invariant we will now
introduce 
Gelfand-Kirillov dimension, which basically measures the growth of a
ring.
Fix a    finitely generated $k$-algebra $R$. Pick generators $r_1,\dots,
r_m$  for $R$ and define an increasing filtration on $R$ by
 $\Gamma_1=k+\sum_1^m r_ik$ and $\Gamma_j=\Gamma_1^j$
for $j>1$.   
Then  the  {\it Gelfand-Kirillov dimension}\label{GKdim-index}
 of $R$ is defined to
be
\begin{align*}
\GKdim(R) \;=\;& \inf \left\{\alpha\in \mathbb R : \dim_k\Gamma_n \leq
n^\alpha
 \text{ for all } n\gg 0\right\}\\
 \;=\;& \limsup \frac{\ln(\dim \Gamma_n)}{\ln(n)}
 \end{align*}

It is routine to check that $\GKdim(R)$ is  independent of the choice of
generators. Also, for a finitely generated commutative algebra $A$, $\GKdim(A)$
reduces  to the degree of the Hilbert-Samuel polynomial as defined, for
example, in  \cite[Chapter~11]{AM1}, and so $\GKdim(A)$ equals the Krull
dimension of $A$. For these and other  basic facts about this dimension,  the
reader is referred to \cite{KL}.  

As  so often happens, many of the basic properties of
Gelfand-Kirillov dimension are much nastier in the noncommutative
universe than
in the commutative one.  For example, it need not be an
integer; indeed there are examples of Bergman that show that $\GKdim(R)$
can
take any  real number  $\alpha\geq 2$ \cite[Theorem~2.9]{KL}.
 Fortunately,  if
$\GKdim(R)\leq
2$ then  $\GKdim(R) $ equals zero, one or two. (See \cite[Theorems~2.5,
2.9
and  Proposition~1.4]{KL}.)  Curiously, though, there are no examples of
domains with non-integer Gelfand-Kirillov dimension. 

If $\GKdim(R)=0$, it is an easy exercise to show that $R$ is a finite
dimensional $k$-algebra. Thus the first interesting case occurs when
$\GKdim(R)=1$.

\begin{prop}\label{prop1.3} \cite{SmWa}
Let $R$ be a domain, with $\GKdim(R)=1$, that is 
finitely generated as a $k$-algebra.
Then  $R$ is  commutative. 
\end{prop}

If one takes a geometric point of view, whereby an
(ungraded) noncommutative
domain of Gelfand-Kirillov dimension $n$ should be regarded as the
coordinate ring of
a noncommutative irreducible 
affine scheme of dimension $n$, then this result says
 that such a scheme is a genuine irreducible curve. 
Similarly, an irreducible
 noncommutative projective scheme of dimension zero, which would
correspond to a graded domain $R$ with $\GKdim(R)=1$, is just a point. 
Almost nothing is known about higher dimensional affine analogues of this
proposition, but the underlying thesis of this paper is that the projective
analogues do seem to exist. 

Proposition~\ref{prop1.3} actually generalizes to any ring 
of Gelfand-Kirillov dimension one, even over fields that are not algebraically
closed. The notation is explained in \S\ref{notation}.

\begin{proposition}\label{prop1.1} \cite{SSW}
Let $R$ be a finitely generated $k$-algebra and assume that
$\GKdim(R)=1$.
Then  $N(R)$ is nilpotent and $\bar{R}=R/N(R)$ 
is a finite module over its centre $Z=Z(\bar{R})$.
Moreover, $Z$ is a finitely generated  semiprime $k$-algebra of 
Krull dimension one.
In particular,  $Z$ and  $\bar{R}$ are noetherian.
\end{proposition}
 
The reason why Proposition~\ref{prop1.3} follows from Proposition~\ref{prop1.1}
is as follows: If $R$ is a domain, then 
$Z=Z(R)$ is also domain, and its field of fractions
$F=\operatorname{Frac}(Z)$ 
has transcendence degree at most one. Since $k$ is algebraically closed, 
Tsen's Theorem \cite[Exercise 7, p.375]{Cohn} says that 
$F$ is the only  division ring, finite dimensional 
as an $F$-module, with centre $F$. Since $RF$ is just such a
division
ring,  $R$ is commutative and we are done.

\section{Noncommutative Projective
Curves I: Domains of Quadratic Growth.}\label{sec4}

Having dismissed the algebras of Gelfand-Kirillov 
dimension one, the
next case of interest is that of algebras of 
Gelfand-Kirillov
dimension two.  Here, we are   concerned with
 the structure of
finitely graded algebras, since these would 
provide the homogeneous coordinate
rings of the (putative) noncommutative projective curves. 

 In this section, we 
describe all   finitely graded domains of Gelfand-Kirillov
dimension $2$, together 
with their projective geometry.  
The case of prime rings over arbitrary fields will be covered in
\S\ref{sec5}. Roughly speaking, these results show that there are
only
two types of noetherian domains $R$ that can occur: one either has
a twisted homogeneous coordinate ring $B(Y,\calL,\sigma)$,
 in the sense of \S\ref{sec2.2},
or one has an analogue of the characteristic zero 
version of the ring $S=k+Uy$
 from \eqref{non-chi}. In either case, $\QGr(R)$ is essentially 
 $\Qch(Y)$, for some curve $Y$. Thus these results really do
justify the intuition that any such domain is the homogeneous coordinate
ring of
a noncommutative projective curve.

Fix a finitely graded domain $R$ with $\GKdim(R)=2$.
As in commutative algebra, a natural first step towards describing 
 $R$ is to determine its ``birational'' structure; that is, to
determine
the graded ring of fractions  $Q=Q(R)$. This is the  ring  obtained   by
inverting the homogeneous elements in $R$ (for the nontrivial
proof that this ring even exists combine 
\cite[Theorem C.I.1.6]{NVO} with \cite[Theorem 4.12]{KL}).
   Then 
\cite[Theorems A.I.5.8 and C.I.1.6]{NVO}  
show that $Q(R)=D[z,z^{-1};\sigma]$, where    $\sigma$ is an
automorphism of the  division ring $D$ and
$D[z,z^{-1};\sigma]$ is the twisted Laurent polynomial ring; thus,
as a $k$-vector space,
 it is isomorphic to the usual Laurent polynomial ring, but
multiplication is twisted by $za=\sigma(a)z$, for $a\in D$. 
However, much more is true:

\begin{thm}\label{thm3.1} \cite[Theorem~0.1]{Staf5}
Assume that  $R$ is a finitely graded domain  with  $\GKdim(R)=2$ and
 write   $Q(R)= D[z,z^{-1};\sigma]$ for the graded
quotient ring of $R$.

\begin{enumerate}
\item[(i)]  $D=K$ is a field that is a finitely generated 
and of transcendence degree  one as an extension of $k$.

\item[(ii)]  Conversely, if $Q=Q(R)$ is as in (i) and if $S$ is any
finitely graded subring of $Q$, then $\GK(S)\leq 2$.
Indeed, there even exists a constant $c$ such that $\dim_kS_n\leq
cn$, for all $n\geq 1$.
\end{enumerate}
\end{thm}

The idea behind the proof is to show that $\GKdim(T)\leq 1$,
 for any finitely
generated subalgebra $T\subseteq D$. 
 Most of the theorem then follows  from
Proposition~\ref{prop1.3}.

 The notation of Theorem~\ref{thm3.1} 
will be fixed for the rest of the 
section. By the theory of algebraic 
curves \cite[\S I.6]{H},
 $K=Q(R)_0$ is the
function field  of a  unique smooth projective
curve $X$ and the automorphism $\sigma$ induces an automorphism,
which will still be denoted by $\sigma$, on $X$. Thus
the theorem  can  really be regarded as saying  that $R$
is determined birationally by an automorphism of that curve. 
 
However we can give a much more detailed
 geometric description of the ring 
$R$ and  the model
for our program is again provided by algebraic geometry, in the form of
 Serre's Theorem~\ref{serre}. 
 Remarkably, this   generalizes completely:

 \begin{thm}\label{thm3.3} \cite[Theorem 0.2]{Staf5}
 Let $R$ be a connected graded domain of GK-dimension
$2$ and assume
 that $R$ is  generated by elements of degree $1$.  Then there
exist:

\begin{enumerate}
\item[(i)] A projective curve $Y$, birational to $X$ and
with an induced action of $\sigma$,

\item[(ii)] An ample invertible sheaf $\calL$ on $Y$,  \end{enumerate}
\noindent such that  $R$ embeds into the twisted homogeneous coordinate 
ring
$B=B(Y,\calL,\sigma)$ with $\dim_k(B/R)<\infty$. 
 Equivalently,  $R_n\cong
\HH^0(Y, \calL\otimes \calL^{\sigma}\otimes\cdots\otimes
\calL^{\sigma^{n-1}})\quad
\text{for } n\gg 0.$ 
\end{thm}

Since ample invertible sheaves on curves are automatically
$\sigma$-ample,  we
may immediately apply Theorem~\ref{AZZ} to obtain:

 \begin{cor}\label{cor3.4} Let $R$ be as in Theorem~\ref{thm3.3}.
Then $R$ is a noetherian domain, and the quotient category $\qgr(R)$ 
is  equivalent to the category $\coh(Y)$ of coherent
sheaves on $Y$.  \end{cor}

Since twisted homogeneous coordinate rings are so well understood, 
one can use the theorem to give further results about the structure of
$R$. 
For example, under the equivalence of the corollary, two-sided ideals of
$R$
correspond to the (ideals of definition of) closed  $\sigma$-invariant 
subschemes of ${{Y}}$. This quickly leads to: 

\begin{cor}\label{cor3.5} Let $R$ be as in Theorem~\ref{thm3.3}. 

\begin{enumerate}
\item[(i)] 
If the order $|\sigma|$ of $\sigma$ is finite, then $R$ is
a finite module over its centre.

\item[(ii)] If $|\sigma|=\infty$, then ${{X}}$ is either rational or 
elliptic. In either case, $R$ is a primitive ring; that is, $R$ has 
a faithful simple right module.  

\item[(iii)]  If ${{X}}$ is rational and $|\sigma|=\infty$, then 
 $R$ has just one or two   
 graded prime ideals $P_i$ of height one  and, in
each
 case, $R/P_i$ is a domain.
 All other non-zero prime ideals $Q$ of 
$R$ are of finite codimension in $R$.

\item[(iv)] If ${{X}}$ is elliptic and $|\sigma|=\infty$,
 then ${{Y}}={{X}}$ and the only non-zero prime ideal of $R$ is the 
 augmentation ideal $R_+=\bigoplus_{i>0}R_i$.
 \end{enumerate}
 \end{cor}

The case that $R$ is not generated in degree one, however, differs
essentially
from commutative algebra, and exhibits a striking new phenomenon.

 To see this, define the
$d$-th Veronese ring\label{veronese-index} of $R$ to be the graded ring
$R^{(d)}=\bigoplus_{n\geq
0}R_{nd}$, graded by $R^{(d)}_n=R_{nd}$. If $R$ is generated in degree
one then it is easy to prove that $\QGr(R)\cong
\QGr(R^{(d)})$. So from a geometrical standpoint $R$ and $R^{(d)}$ are
equivalent.  
Even if $R$ is not generated in degree one then
it may happen, as is always the
case
if $R$ is commutative, that some Veronese ring $R^{(d)}$ is generated in
degree
one. If so, then $R$ is a finite $R^{(d)}$-module and, as in the
commutative
case, the structure of $R$ is closely related to that of $R^{(d)}$. 
Theorem~\ref{thm3.3}  can be applied to such rings and, after the
appropriate
modifications, analogues of Corollaries \ref{cor3.4} and \ref{cor3.5} 
will  hold for $R$. The
details can be found in \cite[\S 6]{Staf5}.

Unlike the commutative case, though, it may also happen that {\it no}
Veronese ring is generated in degree one; indeed,
Theorem~\ref{thm3.7}(i)
 can be regarded as  saying that this is the generic case when
 $|\sigma|=\infty$. An explicit 
 example is given, over any field $k$,
  by the ring $S=k+Uy$ from
\eqref{non-chi} (see \cite[Corollary~3.2]{Staf4}).
The next two theorems from \cite{Staf5} summarize
 the main results about this class of rings and, in particular, explain
the
 dichotomy \eqref{s-is-noether} 
 that $S=k+Uy$  is noetherian if and only if $\charact k = 0$.

\begin{thm}\label{thm3.6} \cite[Theorem 0.4]{Staf5}
Assume that  $R$ is a finitely graded domain of
GK-dimension $2$.

\begin{enumerate}
\item[(i)] If $|\sigma|=\infty$, then $R$ is a noetherian 
primitive ring, and every Veronese ring of $R$ is finitely generated.

\item[(ii)] 
 Assume that  $|\sigma|=d<\infty$ and 
 that no Veronese ring
of $R$ is generated in degree one. Then $R$ is not noetherian.
Indeed, $R^{(d)}$ is a commutative ring which  is not
finitely generated.
\end{enumerate}
\end{thm} 

 Part (ii) of the theorem is easy: Since $\sigma ^d=1$, 
$R^{(d)}$ is a subring of the commutative ring $K[z^d,z^{-d}]$, and so
is
commutative. If $R^{(d)}$ were finitely generated, then classical
results ensure
that some further Veronese ring is generated in degree one. Thus 
the interesting case of the theorem is when $|\sigma|=\infty$.
 In this case, much more can be said about the structure of $R$,
 but for simplicity we only
state the result for a Veronese.    A
complete description of the rings which arise in case (ii) 
of Theorem~\ref{thm3.6} remains to be found.

\begin{thm}\label{thm3.7} \cite[Theorem 0.5]{Staf5}
Let $R$ be as in Theorem~\ref{thm3.6}.  Assume
that
no Veronese ring of $R$ is generated in degree one, and that
$|\sigma|=\infty$. Then, after replacing $R$ by an appropriate Veronese
ring, the following hold:
 
 \begin{enumerate}
\item[(i)]  There is an algebraic curve ${{Y}}$, birational to $X$ and 
   with an
induced action of $\sigma$,
 and two   invertible sheaves $\mathcal L \subsetneq
\mathcal N$ on $Y$, such that
$$R_m=\text{H}^0({{Y}},\,\mathcal L\otimes \mathcal N^\sigma \otimes
\cdots\otimes \mathcal N^{\sigma^{m-1}}), \qquad{\rm for\ all\ } m\geq 1.$$  
Moreover, $\Lscr=\Nscr$ locally on every finite orbit of $\sigma$.
 
\item[(ii)]  $R$ is generated in degrees one and two.

\item[(iii)]  $R$ is contained in the twisted homogeneous coordinate
ring $B=B({{Y}},\,\mathcal N,\sigma)$. Moreover, $B$ is
finitely generated as a right $R$-module but infinitely generated as a
left $R$-module.

\item[(iv)]  $R$ does not satisfy $\chi_1$. However, 
 $\qgr(R)$ is equivalent to  $\coh(Y)$. 
\end{enumerate}
\end{thm} 

The converse to part (i) holds in the sense that, if $\Lscr$ and $\Nscr$
are as
defined there,  then $R=k\oplus \bigoplus_{n\geq 1} 
\text{H}^0({{Y}},\,\mathcal L\otimes \mathcal N^\sigma \otimes
\cdots\otimes \mathcal N^{\sigma^{n-1}})$ 
is a noetherian  connected graded domain.
Moreover, if $\Lscr\not= \Nscr$ then no Veronese ring $R^{(d)}$ is generated in
degree one.

{\it Exercise:}  Write the ring $S=k+Uy$ from \eqref{non-chi} in the form
described by Theorem~\ref{thm3.7}.
Under the appropriate identifications, you should find that
$\Lscr=\Oscr_{\PP^1} \subset \Nscr=\Oscr_{\PP^1}(1)$ and that
$B(Y,\Nscr,\sigma)=yUy^{-1}\cong U$.

One curious consequence of the  results 
from this section is that they
contradict the pervasive feeling that 
noncommutative rings are more varied and
complicated than their commutative 
counterparts; notably that one should 
expect that the rings that are farthest from
commutative should have the worst properties.
For connected graded domains of dimension
two, the opposite is  true: The only nasty (non-noetherian) examples are close
to commutative  (formally, they embed into a finite free
module over a commutative subalgebra),
 while any  such domain that  is not
close to commutative is defined in 
terms of a curve of genus $\leq 1$.

Let us at least indicate how to construct the curve $Y$ in
Theorems~\ref{thm3.3} and \ref{thm3.7}.    We will only do
this
in the case when
$R$ is generated by $R_1$, but the general case is similar.
Fix a non-zero element  $z\in R_1$ and write  $R_m=\RR_mz^m$, so that
 $\RR_m$ is a
subspace of the function field $Q(R)_0=k(X)$. Multiplication in $R$ 
translates to the statement that 
\begin{equation}\label{equ3.8}
\RR_n\sigma^n(\RR_m)\subseteq \RR_{n+m}\qquad\text{for } n,m\geq 0. 
\end{equation} 
The space $\RR_m$ defines a
morphism $\pi_m: X\to   \mathbb P(\RR_m^*)$  
\cite[Theorem~II.7.1]{H} and we set $Y_m=\Im(\pi_m)$.
A key step in the proof is to show, possibly after replacing $R$ by a
Veronese ring $R^{(d)}$, 
 that $Y_m\cong Y_{m+1}$ for all   $m\gg 0$. Once this has been
 achieved, then one simply takes $Y=Y_m$, for some large $m$.
Set $\calL_m= \O_{Y}\RR_m$, which  is an  invertible sheaf on $Y$.
A second key fact is that
 $$\RR_m=\text{H}^0(Y,\, \calL_m)\qquad\text{ for all }m\gg 0.$$  
We are now done. Indeed, since $R$ is generated by $R_1$, \eqref{equ3.8}
implies that
$\RR_m=\RR_1\RR_1^\sigma\cdots \RR_1^{\sigma^{m-1}}$ and so $\mathcal
L_m=\mathcal L_1\otimes\cdots\otimes\mathcal L_1^{\sigma^{m-1}}$. 
Thus  Theorem~\ref{thm3.3} 
holds with $\calL=\calL_1$.

\section{Noncommutative Projective Curves II: Prime Rings.}\label{sec5}

In this section, we outline how to generalize the results of the last
section to prime rings of quadratic growth. 
There are two reasons for this. 
Despite the fact that most of the results 
in this survey are concerned with domains, the 
correct context is probably that of
prime rings---certainly this is true for the most significant 
applications of ring theory---and so this section can be thought of as 
illustrating how such generalizations might go. 
Secondly, one of the best understood areas of 
noncommutative algebra is that of 
classical orders (which we  define to be
 non-graded prime rings of
Gelfand-Kirillov dimension one). Given the success of
\S\ref{sec4}, one should expect 
a similar generalization of the theory of
classical orders to the projective setting and  this 
section provides just such a
theory. Once  one works with prime rings rather than domains, 
the standing assumption that $k$ be algebraically closed
provides no advantages and so  the results of this section actually hold for
arbitrary base fields.  All the results of
this section are taken from \cite{Staf6}.

Thus, we wish to study  finitely graded prime rings $R$ with $\GK(R)=2$.
Although the results of the last section generalize naturally, there are
some differences and to explain this
 we begin with a series of examples. 
In all these examples, 
 we begin with the standard
quantum plane  $U=\mathbb C\{ x,y\}/(xy-yx-x^2)\cong B(\PP^1, \Oscr(1), \sigma)$,
as in \eqref{babyeg}, but we now  consider 
subrings of the
$2\times 2$ matrix ring $V=M_2(U)$. 

The first obvious difference is that we have to use, as our basic building
block,  rings $B(\Ga, \cal B_1, \tau)$ defined 
 in terms of orders $\Ga$ over a curve $Y$ rather
than the twisted homogeneous coordinate rings $B(Y,\Lscr,\sigma)$ that were
used in \S\ref{sec4}.  This does not cause any serious complications;
in fact it is even an advantage since it allows for more subtle examples from
the same basic construction.
The construction is given in detail in \cite{Staf6}, but 
the idea is as follows. Assume that $X$ is a  projective variety 
and that $\Escr$ is an order inside a central
semisimple $k(X)$-algebra $A$. Let
 $\tau$ be a $k$-algebra
automorphism of $A$ that restricts to an automorphism
$\sigma$ of $k(X)$ and $X$ and assume that 
 $\Lscr$ is an invertible sheaf of $(\cal E,\cal
E^\tau$)-bimodules   contained in $A$. 
Modulo these definitions, 
the {\it twisting ring}\label{twistingring-index}
  $B =B(\Escr,\Lscr,\tau)$ looks just 
like the twisted homogeneous coordinate ring
 $B(X,\Lscr,\sigma)$.
Indeed,   
$$B \cong \bigoplus_n
\text{H}^0({{X}}, \calL_n) \qquad{\rm  where} \quad \calL_n
= \calL\otimes_{\mathcal E^\tau} \calL^\tau 
\otimes_{\mathcal E^{\tau^2}}\cdots \otimes_{\mathcal
E^{\tau^{n-1}}} \calL^{\tau^{n-1}},$$ and  multiplication is
defined by the analogue of
\eqref{equ2.2}. If $\Lscr$ is ample as an $\Oscr_X$-module, then 
 $s=\tau_\ast(-\otimes_{\Escr} \Lscr)$ is  ample in the sense of
 Definition~\ref{ampleness}, and so
Theorem~\ref{AZ} can again be used to describe $B$.
 
 As an exercise, the reader may like to prove that the matrix ring $V=M_2(U)$
 satisfies $V\cong
 B(\Ga, \cal B_1,\tau)$, where $\Ga=M_2(\cal O_{X})$, $\cal
 B_1=\cal O_{X}(1)\otimes_{\cal O_{X}}\Ga$ and $\tau $ 
 is the natural extension to $\Ga$ of the automorphism 
 $\sigma$ from~\eqref{babyeg}.
 A more subtle example appears at the end of the section.

A more serious complication occurs if one tries to generalize
Theorem~\ref{thm3.3} to the prime case, since it is now easy to
construct rings
 generated in degree one that are not noetherian, let alone twisting
rings.
Here are a couple of simple examples.

\begin{example}\label{prime-a} Consider 
\begin{equation*}\label{equ4.5}
R'=\begin{pmatrix} U&\hfill Ux \\U& \mathbb C +Ux \end{pmatrix}. \end{equation*}
Then it is not hard to show \cite{Staf4}
 that $R'$ is neither left nor right noetherian.
However, it is generated in degree one.
In fact, this example still works if we replace $U$ by the commutative
polynomial ring $\mathbb C[x,y]$. 
\end{example}

\begin{example}\label{prime-b}  Consider 
\begin{equation*}\label{equ4.55}
R''=\begin{pmatrix} U&\hfill Uy \\U& \mathbb C +Uy \end{pmatrix}.
 \end{equation*}
Then $R''$ is both left and right noetherian and 
  it is still generated in degree one (see \cite{Staf4}, again).
However, it follows from \eqref{equ4.6} below that $R''$ is not a
twisting ring.   
\end{example}

Thus generation in degree one is not very helpful and so one needs a
more
subtle geometric  conditions to determine the twisting
rings, or even when a ring is noetherian. 
These are provided by \eqref{equ4.6} and Theorem~\ref{thm4.7}.  
Observe that  $R'=\mathbb C+I'$, 
and $R''=\mathbb C+I''$,
 for  left ideals $I'$ and $ I''$ of $V$. 
The difference between them is that $V/I'$ is supported (in a manner
that  is
 made precise in Theorem~\ref{thm4.1})
 in $\PP^1$ at the fixed point of
 $\sigma$, whereas $V/I''$ is supported on an infinite
orbit
 of $\sigma$.
  Thus  $R''$ looks very like---and behaves very
 like---the ring $\mathbb C+Uy$ that played such a  
 prominent r\^ole in the last   
section. In contrast, the analogue to $R'$ would be $\mathbb C+Ux$, 
which is not
finitely generated and so never appeared in 
\S\ref{sec4}.

With these examples to hand, the main result of this section is easy to
state
informally: {\it A finitely graded  prime ring $R$, with $\GKdim(R)=2$,
is
noetherian if and only if some Veronese ring 
 $S=R^{(d)}$ is either a twisting ring $B=B(\Ga, \cal
B_1,\tau)$ or else of the form $S=k+I$, where $I$ is a left ideal of $B$
such
that $B/I$ is supported on a union of infinite orbits of $\sigma$.
In either case, $\qgr(S)\simeq \coh(\Ga)$.}

We now formalize these ideas. It  is actually easier
to work with semiprime rather than prime rings, so fix  a  base field $k$
and let $S=\bigoplus_{i\geq 0} S_i$ be a   semiprime  finitely
graded
$k$-algebra with quadratic growth. 
Unlike the last section, the assertion that the quotient ring $Q(S)$ exists 
is a genuine assumption and is one of the equivalent definitions that $S$ be a
Goldie ring\label{goldie-index} (see 
\cite[Theorem~C.I.1.6]{NVO} for the details and note that semiprime
noetherian rings are automatically Goldie).
  For simplicity, we will
make two further assumptions: Define  $S$ to be {\it nice}\label{nice-index} if  
\begin{enumerate} \item[(a)]  $\GK(S/P)=2$ for
all minimal prime ideals $P$ of $S$ and 
\item[(b)]  $S$ contains a non-zero divisor
$z\in S_1$. 
\end{enumerate}These assumptions are only  being made to simplify our 
notation,
and can be easily removed (see \cite{Staf6}). Indeed, (a)  will
always hold if $S$ is prime, while (b) can always  be assumed by
passing  to a
Veronese ring.

One advantage of condition (b) is that the graded quotient ring 
$Q=Q(S)$ of
$S$ can now be written  $Q = A[z,z^{-1};\tau]$, 
where the semisimple Artinian algebra
$A$ equals $Q_0$ 
and $\tau $ is an automorphism of
$A$ \cite[Lemma~A.II.1.5]{NVO}.  Standard techniques  from Goldie's theory 
(see \cite[\S 1]{AS})
allow one to generalize  Theorem~\ref{thm3.1} 
to this context and they show
that {\it $A$ is a finite module over 
its centre $K$ and $K$ is a direct sum
of fields of transcendence degree one.}  
Let $X$ denote the smooth projective 
 model of $K/k$
and write $\sigma$ for the restriction of $\tau$ to $K$ and $X$.  
We fix this notation throughout this section.

Pick a $\sigma$-stable projective model ${{Y}}$ of $K/k$, an
$\O_{{{Y}}}$-order
$\Ga$ in $A$ and an ample invertible $(\Ga, \Ga^\tau )$-bimodule
$\mathcal
B_1$.  The next result provides the general version of
Example~\ref{prime-b},
and can also be thought of as the natural generalization of the
construction of
Theorem~\ref{thm3.7}(i).  

 \begin{thm}\label{thm4.1} Let ${{Y}}$, $\Ga$ and $\mathcal B_1$ be as
above.
   Let $\mathcal R_1$ be an essential right $\Ga^\tau$-submodule of
$\mathcal B_1$ such that $1\in \HH^0({{Y}}, \mathcal R_1)$ 
and such that $\mathcal
R_1=\mathcal B_1$ locally on every finite orbit of $\sigma$ on ${{Y}}$. 
For
$n\geq 1$, set
$$\mathcal R_n =\mathcal R_1\otimes_{{\Ga}^\tau}\mathcal
B_1^\tau\otimes_{{\Ga}^{\tau^2}}\cdots\otimes_{\Ga^{\tau^{n-1}}} \mathcal
B_1^{\tau^{(n-1)}}.$$ 
   Then ${R} = k+\bigoplus_{n\geq 1} \HH^0({{Y}},\mathcal R_n)$ is
a noetherian \order{}.  
\end{thm}

{\it Exercise:} Write the   ring $R''$ from Example~\ref{prime-b}
 in this form. The answer, modulo a left-right switch, can be found in 
 \cite[Example~5.13]{Staf6}.   

Once again, the rings constructed by Theorem~\ref{thm4.1} have a
pleasant
module structure.
\begin{cor}\label{cor4.2} 
Let $R$ be the ring constructed by 
  Theorem~\ref{thm4.1}. Then
 \begin{enumerate}
\item[(i)] 
 The categories $\qgr(R)$  and $\coh(\Escr)$ are equivalent.
\item[(ii)] 
   $R$ is a 
subidealizer in $B=B(\Ga, \mathcal B_1, \tau)$ in the sense that,
for some $r\geq 1$,  
 $V=\bigoplus_{i\geq r}R_i$ is a right ideal of $B$. 
Moreover, $B$ will be a finitely generated right $R$-module,
although it may be infinitely generated as a left
$R$-module.  
\end{enumerate}
\end{cor}

The main result of \cite{Staf6} is the following partial converse to
Theorem~\ref{thm4.1}.
 
\begin{thm}\label{thm4.3} Let $S$ be a noetherian \order{}. For
some $n$, the Veronese ring $R=S^{(n)}$  has the form
described in Theorem~\ref{thm4.1}.
Moreover, if $\mathcal S_n = \mathcal
O_{{Y}} S_n$, then $S_n = \HH^0({{Y}},\mathcal S_n)$ for all $n\gg 0$. 
 \end{thm}

If $S$ is a noetherian \order,  as in Theorem~\ref{thm4.3}, then the
module
 structure of $S$ may be more complex than that described in
 Corollary~\ref{cor4.2}---but then  this is also true for a commutative
ring 
 that is not generated in degree 
 one. Nevertheless, the equivalence 
nearly extends to this case.  Formally,  by 
\cite[Proposition~6.1(iv)]{Staf5}
 the category of 
 graded  $1$-critical $S$-modules that are generated in degrees 
 $n\mathbb Z$ is equivalent to the category of $1$-critical 
 modules over $R=S^{(n)}$ and so  one may apply Corollary~\ref{cor4.2}
 via this correspondence. Consequently, the 
 simple objects in $\qgr(S)$ that are generated in degrees 
 $n\mathbb Z$ are in one-to-one correspondence with the
 simple $\Ga$-modules.

 The final question we need to address is to determine geometrically
 when a \order{} is noetherian.  
  To describe this condition, it
is convenient to write a \order{} $R$ as $R=\bigoplus \RR_iz^i$,
so that $\RR_i \subset A$ for all $i$.  The fact that $R$ is a graded
ring translates to the inclusions $\RR_i\RR_j^{\tau^i}\subseteq
\RR_{i+j}$.  Let $N$ denote the reduced norm map
$A \to K$.  For a point $p\in X$ and a subset $V$ of $A$, we
denote by $\rr(V;p)$ the maximum order of a pole at $p$ of $N(v)$, for
$v\in V$. The key condition is
\begin{equation}\label{equ4.6}
\rr(\RR_i;p)+\rr(\RR^{\tau^i}_j;p)
=\rr(\RR_{i+j};p). \end{equation}

If $A=K$, then \eqref{equ4.6} will hold whenever $\RR_i\RR_j^{\tau^i}
=\RR_{i+j}$. However, this is  not true in  
general---which should come as no surprise, since it is related to the fact that 
 the determinant of a sum (of products) of
  matrices cannot be described in  terms of the
determinants of the individual matrices. 
In particular, it fails  for the rings $R'$ and $R''$ of
Examples~\ref{prime-a} and \ref{prime-b}.    If one computes
$\rr(\RR_i;p)$ for  $R'$ and $R''$, one finds that \eqref{equ4.6} fails for
$R'$ at just one point, the  fixed point of $\sigma$, whereas for $R''$ it
  fails only on an infinite orbit of $\sigma$. The details can be found in
\cite[Example~5.14]{Staf6}. Thus  the following result should come as no
surprise.

\begin{thm}\label{thm4.7} Let $S$ be a  \order. Then:

\emph{(i)}   $S$ is right noetherian if and only if it is left
noetherian. This
is true if and only if  there is an integer $n$ such that   $R=S^{(n)}$
is a
\order{} for which  \eqref{equ4.6} holds on all finite orbits of
$\sigma$ on
$X$.

\emph{(ii)} Some Veronese ring $S^{(d)}$ is a twisting ring  $B(\Ga,
\cal
B_1,\tau)$  if and only if some   $R=S^{(n)}$ is a \order{} for which 
\eqref{equ4.6} holds at all $p\in X$.\end{thm}
 
Many of the standard examples in PI theory can be modified 
to be graded rings of Gelfand-Kirillov dimension two, and so 
can be reinterpreted in terms of twisting rings. As an example, we invite the
reader to describe the following algebra as a twisting ring:
\begin{equation} \label{equ2.4} 
 B\; =\; \left\{ \begin{pmatrix} \hfill a
+yb_{11}&  \hfill  x{\frac{ \partial(a)}{\partial x}} + yb_{12} \\  \hfill 
yb_{21} &\hfill a+yb_{22}\end{pmatrix} :  a, b_{ij}\in k[x,y]\right\}. 
\end{equation} 
One solution is $B\cong B(\mathcal E, \mathcal B_1, \tau)$,
where  $\mathcal E= \mathcal O_{\PP^1}  + M_2(\mathcal O(-1))$ 
and
$\mathcal B_1 = \mathcal  E\otimes_{\mathcal O}\mathcal O(1)$.  
One interesting feature of this example is that $\tau$ is the inner
automorphism given by conjugation by  $\left(\smallmatrix
1&1\\0&1\endsmallmatrix\right)$.
 The details can be found in \cite[Example 5.14]{Staf6}.
   This  is one of the standard examples of a noetherian  ring
that satisfies a polynomial identity but  is not a finite module over its
centre  \cite[Example~5.1.16]{Rw}.   
 
\section{Noncommutative Smooth Proper Curves}
\label{secy1} 

Although the results from the last two sections are very
satisfying they should be only part of the story. As in the commutative
case one would very much prefer to define noncommutative curves
abstractly and then to show that they can indeed be described by
graded rings.  We provide one possible method in 
this section.

This has three steps: first, we describe the 
 homological 
properties that are satisfied by a  smooth proper curve
 and which should be satisfied by 
any noncommutative analogue of such a curve. 
Secondly, we  provide a list of abelian 
 categories that  satisfy these conditions.  
Finally, we show that these categories are 
indeed determined by those
properties.

\subsection{Geometric Conditions on Abelian Categories}
\label{sec6}

In this subsection, we assume that $X$ is 
\emph{smooth, proper and connected} $k$-scheme of dimension $n$  
and we isolate the
homological conditions that we want to be satisfied by our 
noncommutative analogues of
$\Cscr=\coh(X)$. For a discussion of the singular case we refer to
\S\ref{secmorecomments}. 
  One could also have associated  $\Qch(X)$
to $X$ but,  as Theorem~\ref{Gabriel} shows, 
this would not make too much difference.

  The category
$\Cscr=\coh(X)$ clearly satisfies the following conditions.
\begin{itemize}
\item[\bf(C1)] $\Cscr$ is noetherian.
\item[\bf(C2)] $\Cscr$ is {\it $\Ext$-finite}\label{extfin-index}. 
In other words,
the Yoneda $\Ext$ group  satisfies $\Ext^i(A,B)<\infty$
  for all $A,B\in\Cscr$ and   all $i$.
\item[\bf(C3)] $\Cscr$ has {\it homological dimension $n$}\label{homdim-index}
 in the sense
that
  $\Ext^i(A,B)=0$ for $A,B\in\Cscr$ and for $i>n$, and $n$ is minimal
  with this property.
\end{itemize}
Unfortunately, far too many noncommutative 
rings have finite 
homological dimension and so these three conditions are not nearly
sufficient.
A fundamental, but slightly more subtle property 
of smooth proper schemes 
is the Serre duality theorem. This asserts that there exists a dualizing
sheaf $\omega$
such that  for $\Fscr\in\Cscr $
there are natural isomorphisms
\[
H^i(X,\Fscr)\cong \Ext^{n-i}(\Fscr,\omega)^\ast
\]
A very elegant reformulation of Serre duality was given  
 in \cite{Bondal4}. Let $D^b(\Gscr)$ denote the derived category of
bounded complexes over an abelian category 
$\Gscr$. A \emph{Serre functor}\label{serre-fn} on $D^b(\Gscr)$ is 
an autoequivalence
  $F:D^b(\Gscr)\r D^b(\Gscr)$ such that there are isomorphisms
\begin{equation}
\label{serreisos}
\Hom(A,B)\cong \Hom(B,FA)^\ast
\end{equation}
which are natural in $A,B\in D^b(\Gscr)$. Then 
 Bondal and Kapranov show
that 
Serre duality can be reinterpreted as saying:
\begin{itemize}
\item[\bf(C4)]  $\Cscr$ satisfies \emph{Serre duality}
 in the sense that there exists a
  \emph{Serre functor} on $D^b(\Cscr).$
\end{itemize}

\subsection{Some examples}
\label{secy}
An obvious test question for noncommutative geometry would be to 
to classify abelian categories  satisfying  the homological properties
(C1)--(C4) listed above. Since this is far too optimistic a project at
present, we just consider the 
{\it hereditary} categories\label{heredcat-index}; that is the
ones 
of homological dimension one. Here, at least one has a chance, since
there is a
small list of examples. 

 The most obvious   \emph{hereditary noetherian  $k$-linear
$\Ext$-finite
 abelian categories satisfying Serre duality}  are:
  
 \begin{itemize}
\item[\bf(E1)] Let $X$ be a smooth projective connected curve over $k$
  with function field $K$ and let $\Oscr$ be a sheaf of hereditary
  $\Oscr_X$ orders in $M_n(K)$ (see \cite{reiner}). Then one proves
  exactly as in the commutative case that $\Cscr_{1}=\coh(\Oscr)$
satisfies
  Serre duality. The Serre functor  is given by tensoring with
  $\uHom(\Oscr,\omega_X)[1]$. 
\end{itemize}

 Three rather special examples, for which it is easy to check 
(C1)--(C4) are:
\begin{itemize}
\item[\bf(E2)] 
$\Cscr_{2}$ consists of the finite dimensional nilpotent
  representations of the quiver  
$\tilde{A}_n$ \cite{Ringel4} (possibly with
  $n=\infty$) with all arrows oriented in the same direction. The
  Serre functor is
  given by rotating one place in the direction of the arrows and
  shifting one place to the left in the derived category.
\item[\bf(E3)]
 $\Cscr_{3}
 =\gr_{\ZZ^2}(k[x,y])/(\text{tor})$ consists of the 
finitely generated 
 $\ZZ^2$-graded $ k[x,y]$-modules modulo 
the finite dimensional modules.  The Serre functor is given by the shift $M\mapsto
M(-1,-1)[1]$ 
\item[\bf(E4)]  The
category $\Cscr_{3}$ has a natural automorphism $\sigma$ of order two
which sends a graded module $(M_{ij})_{ij}$ to $(M_{ji})_{ij}$ and
which exchanges the $x$ and $y$ action. Let $\Cscr_{4}$ be the category
of
$\ZZ/2\ZZ$ equivariant objects in $\Cscr$, i.e. pairs $(M,\phi)$ where
$\phi$ is an isomorphism $M\r \sigma(M)$ satisfying $\sigma(\phi)\phi=
\Id_M$. Properties (C1)--(C4) follow easily from the fact that they hold
for
$\Cscr_{3}$.  In particular the Serre functor is obtained from the one
on $\Cscr_3$ in the obvious way. As described, this construction requires 
$\text{char}\, k\not=2$, but this can be circumvented \cite{ReVdB1}.
\end{itemize}

 We now come to the most subtle example. 
\begin{itemize}
\item[\bf(E5)]  Let $Q$ be a connected  locally
  finite quiver such that the opposite quiver has no infinite oriented
  paths.
    For a vertex $x\in Q$ we have a corresponding projective
  representation $P_x$ and an injective representation $I_x$ and by
  our hypotheses $I_x$ is finite dimensional.  Let $\rep(Q)$ be the
  finitely presented representations of $Q$. 
  The category of all representations of $Q$ is hereditary
and hence
  $\rep(Q)$ is an abelian category. It is easy to check 
   that $I_x\in\rep(Q)$. Hence the functor $P_x\mapsto
  I_x$ may be derived to yield an endofunctor
  $F:D^b(\rep(Q))\r D^b(\rep(Q))$. This functor $F$ behaves like a
  Serre functor  in the sense that we have natural isomorphisms as
  in \eqref{serreisos}, but unfortunately $F$ need not be an
autoequivalence.
  However, there is  a formal
  procedure to invert $F$ so as to obtain a true Serre functor
  \cite{ReVdB1}. This yields a hereditary category $\wrep(Q)$ which
  satisfies Serre duality. Under a technical 
 additional hypothesis (see
  \cite{ReVdB1}) $\Cscr_{5}=\wrep(Q)$ turns out to be noetherian.
\end{itemize}

 Recall that an abelian category
$\Gscr$ is connected if it cannot 
be non-trivially written as a direct
sum $\Gscr_1\oplus\Gscr_2$. The examples (E1)--(E5)
 are all connected and,
remarkably, they are the only ones:  

\begin{theorems}
\label{theoremb} {\em(}\cite{ReVdB1}{\em)}
Let $\Cscr$ be a connected  noetherian $\Ext$-finite
  hereditary category satisfying Serre duality over
   $k$.  Then  $\Cscr$ is one
  of the categories (E1)-(E5).
\end{theorems}

Space does not permit us to sketch the fairly involved proof of this
 theorem. Suffice it to say that, after a series of reductions
which in particular eliminate cases (E2)-(E5), 
 one arrives at the situation where  Theorem~\ref{AZ} can be invoked
to show  that $\Cscr=\qgr(R)$ for some  graded
domain $R$ of quadratic growth. Then, amongst other facts,
Theorem~\ref{thm3.1}
is used to show that $\Cscr$ is of the form (E1).

All the hypotheses for Theorem \ref{theoremb} are necessary. For
example, the non\-commuta\-tive curves considered in \cite{SmithZhang} 
(see  \S\ref{sec7} below)
satisfy (C1)--(C3) but do not satisfy Serre duality.
Nevertheless, it is tempting to ask whether a result similar to 
Theorem~\ref{theoremb} remains valid without the noetherian, or perhaps
the  Serre duality hypotheses. This may require working  up to
derived equivalence.

It is probably too early to be sure what are the ``right'' hypotheses
to demand of a noncommutative analogue of a smooth proper variety.
However, one indication that (C1)--(C4) may need to be supplemented is
that, for example, they do not distinguish between algebraic and
analytic (smooth  compact) surfaces.\footnote{The reason for the 
restriction to surfaces
  is that we do not know if for analytic varieties
the Yoneda
  $\Ext$-groups in $\coh(X)$ coincide with the traditional
  $\Ext$-groups in the category of all $\Oscr_X$-modules. Therefore we
  do not know if Serre duality in the sense of \eqref{serreisos} holds 
for $\coh(X)$. For smooth
  compact analytic surfaces the equality of these
  $\Ext$-groups can
 be proved directly \cite{BondalVdb}.}
   A way
round this problem is to use another subtle property of smooth proper
schemes that was 
also discovered by Bondal and Kapranov, and is not satisfied by
analytic surfaces.

 Recall that a cohomological functor
$H:D^b(\Cscr)\r \mod(k)$ is 
of \emph{finite type} if, for $A\in D^b(\Cscr)$, only a finite number of
the
$H(A[n])$  are non-zero. The appropriate condition is:
\begin{itemize} 
\item[(C5)] Let $\Cscr$ be an $\Ext$-finite abelian category of finite
  homological dimension. Then  $\Cscr$ is
   \emph{saturated}\label{saturated-index} if  every
  cohomological functor $H:D^b(\Cscr)\r \mod(k)$ of finite type is of
  the form $\Hom(A,-)$ (i.e. $H$ is representable).
\end{itemize}
It was shown in \cite{Bondal4} that $\coh(X)$ is saturated when 
 $X$ is a smooth projective scheme,   and that saturation also hold for
categories of
the form $\mod(\Lambda)$ with $\Lambda$ a finite dimensional
algebra. It is easy to show that a saturated category satisfies Serre
duality.

Combined with Theorem~\ref{theoremb}, this gives a much more compact
classification.

\begin{corollarys}
\label{corsat} Assume that $\Cscr$ is a saturated connected
  noetherian $\Ext$-finite  hereditary category.
   Then  $\Cscr$ has one of the following forms:
\begin{enumerate}
\item[(1)] $\mod(\Lambda)$ where $\Lambda$ is an indecomposable 
 finite dimensional  hereditary algebra (this is very special case of
(E5)).
\item[(2)] $\coh(\Oscr)$ where $\Oscr$ is a sheaf of hereditary
$\Oscr_X$
  orders (see (E1) above)   over a smooth connected
  projective curve $X$. 
\end{enumerate}
\end{corollarys}
It is easy to show that the abelian
categories occurring in this corollary  are of the form
$\qgr(R)$ for $R$ a graded ring of Gelfand-Kirillov dimension $\le 2$,
and so 
this result can be regarded as a partial converse to the results of the
last two
sections.  
In a  sense, we have therefore 
generalized the classical commutative result that
smooth proper curves are projective.

Let us briefly turn to more general graded algebras. Although the
saturation
condition seems quite subtle, it is closely connected to the $\chi_n$
conditions we discussed at the end of \S\ref{sec2.1}.
The following result connects the saturation condition with
graded rings.  In this result we define a connected graded ring $R$ to 
have {\it finite right cohomological dimension}\label{cohomdim}
 if the functor 
$\tau M=\dirlim \Hom_{R}(R/R_{\ge n},M)$ 
has finite cohomological dimension. If
 $(\Cscr,\Oscr)=(\qgr(R),\pi R)$, this is the same as demanding that the
 functor $\Gamma(-) = \Hom(\Oscr,-)$ has finite 
 cohomological dimension.

\begin{theorems}\cite{BondalVdb}\label{saturationresult} 
  Let $R$ be a connected graded noetherian ring satisfying
the following hypotheses:  
\begin{enumerate}
\item
  $R$ and its opposite ring $R^\circ$ satisfy $\chi_n$ for all $n$;
\item
$R$  and $R^\circ$ have finite right cohomological dimension;
\item
$\qgr(R)$ has finite homological dimension. 
\end{enumerate}
Then $\qgr(R)$ is
  saturated.
\end{theorems}
We will discuss the first two hypotheses in more detail
at the end of  \S\ref{secmorecomments}. 
The subtlety of the saturation condition is indicated by the following
example.  
\begin{propositions} \label{bv-example}
\cite{BondalVdb} Assume that $X$ is an analytic  K3
  surface with no curves. Then $\coh(X)$ is not saturated.
\end{propositions}

\subsection{More curves}
\label{sec7}
Yet another possible approach to noncommutative 
curves has been given by 
 Smith and Zhang in \cite{SmithZhang}.

If $X$ is a variety over $k$ then the
(closed) points of $X$ are in one-one correspondence with the simple
objects  in $\Qch(X)$. Similarly if $Y\subset
X$ is an irreducible  curve  then its structure sheaf $\Oscr_Y$
is a {\it $1$-critical object}\label{1critob-index} in $\Qch(X)$;
 that is, an object which is not
of finite length but for which every 
quotient object has finite
length (a finite length object\label{finlengthob-index}
 is one that is a finite extension
of simples). 
Now suppose that $\Cscr$ 
is a locally noetherian  Groth\-en\-dieck 
category. 
If we want to view $\Cscr$ as a noncommutative 
variety then, by analogy, we  may argue that the points 
of $\Cscr$ should be the simple
objects in $\Cscr$ and the irreducible  curves in 
$\Cscr$ should be  given by the $1$-critical
objects. Two $1$-critical objects should correspond to the same curve
if they have a common non-trivial subobject. 

Let $L$ be an irreducible curve in $\Cscr$. 
To be consistent with our general
philosophy we should associate an abelian category $\Mod(L)$ 
 to $L$. Following  \cite{SmithZhang},  $\Mod(L)$ is defined to 
  be the smallest
abelian subcategory, containing $L$, of $\Cscr$ which satisfies the
following conditions.
\begin{itemize}
\item $\Mod(L)$ is closed under subquotients and direct sums.
\item$\Mod(L)$ is closed under essential extensions 
by finite length objects. 
\end{itemize}
In \cite{SmithZhang} it is shown that, 
if $(\Cscr,L)$ satisfy suitable
hypotheses which are  special   but 
which nevertheless hold in
many examples, then $\Mod(L)$ is isomorphic to the quotient
 category $\Gscr$ of $\ZZ^t$-graded 
$k[x_1,\dots, x_t]$-modules, modulo those of 
Krull dimension $\leq (t-2)$.
Categories like $\Gscr$ 
 have properties which are very close to
those of $\PP^1$ but they do not satisfy Serre duality when
$t \not= 2$. 
This shows that the Serre duality condition may not always be
satisfied, even in natural situations.

\subsection{Serre duality}
\label{secmorecomments}
It is an unfortunate fact that the 
Serre duality condition (C4) from
\S\ref{sec6} only makes sense 
for categories which have finite
homological dimension, whereas in classical algebraic 
geometry much of the
significance of Serre duality comes in applications to 
singular varieties. We expect the  same to be true 
for noncommutative geometry
and so in this section we describe  a second 
notion of Serre duality that works
more generally. This will naturally lead us 
into a discussion of various homological conditions, 
notably Gorenstein conditions, that will be 
important later.

Suppose that $\Cscr$ is an $\Ext$-finite noetherian  abelian
category with a distinguished object $\Oscr\in\Cscr$.
 By analogy with
the commutative case we write $\Gamma(-)$ for $\Hom(\Oscr,-)$. Then
 $(\Cscr,\Oscr)$ is defined to 
 {\it satisfy classical Serre duality}\label{serredual-index}
  if there exists
an object $\omega^\circ\in D^b(\Cscr)$ together with a natural
isomorphism
\[
\RHom(-,\omega^\circ)\cong {\rm R}\Gamma(-)^\ast
\]
where $\ast$ refers to the $k$-dual. 
  Note that these derived functors can be
  computed using the closure $\widetilde{\Cscr}$, 
 as defined in \S\ref{notation}. 
 Since
$\omega^\circ$ represents a functor   it is clear that it must be
  unique if it exists, 
and, moreover, its existence forces
 $\Gamma(-)$ to have finite cohomological
dimension. 
   Given that   Serre duality is such a powerful tool in the commutative
case, 
 it would be interesting to know how much  more generally it holds. The next
results 
   shows that $\omega^\circ$ does exist for many noncommutative rings.

\begin{theorems}\label{YZ} \cite{YZ}
Let $(\Cscr,\Oscr)=(\qgr(R), \pi R)$ for a noetherian  connected 
graded ring $R$.  Assume  that $R$ satisfies: 
\begin{itemize} 
\item[(a)] $R$ and its opposite ring $R^\circ$ satisfy $\chi_n$ for all
  $n$;
\item[(b)]  $R$ and   $R^{op}$ have finite right cohomological dimension. 
\end{itemize} 
Then $(\qgr(R),\pi R)$ satisfies classical Serre duality.
\end{theorems}

It should be noted that, by invoking the Brown representability
theorem \cite{Keller1,Neeman1}, 
Jorgensen has shown \cite{jorg} that if we allow
$\omega^\circ$ to lie in $D(\QGr(R))$ then it exists under 
far weaker hypotheses. Notably, his result does not require the $\chi$
assumption.

 Jorgensen, Yekutieli, Zhang and others have shown that, under the 
hypotheses of Theorem~\ref{YZ},
 a surprising number of results from the homological
theory of commutative varieties do generalize to the noncommutative
universe.
See, for example, \cite{Jorg2,Jorg3,JZ,StafZ,YZ1,ZW}  and the references
therein.

It is not known how severe are the assumptions of Theorems~\ref{YZ}
and \ref{saturationresult}.   In
particular, it is not known whether every noetherian connected graded ring $R$
has finite cohomological dimension, as would be true if $R$ were
 commutative 
\cite[Theorem III.2.7]{H}. In contrast,
  Theorem~\ref{thm3.7} produces many rings $R$  that do not satisfy the
$\chi$ conditions. One can easily circumvent those examples by
demanding that the noetherian prime ring  $R$ is a 
{\it maximal order}\label{maxorder-index}, in the
sense  that there exists no ring $R\subsetneq R'\subset Q(R)$ with
$aR'b\subseteq R$, for some units $a,b\in Q(R)$. (This is the appropriate
noncommutative analogue of an integrally closed commutative  noetherian
ring.)   Notice that the way we  constructed the example $S=k+Uy$ in
\eqref{non-chi}   ensures that it cannot be a maximal order, and the same is
true of the rings created by  Theorem~\ref{thm3.7}. It remains an open
question whether  every noetherian  connected graded maximal order
 satisfies~$\chi$. 

Here is one case where the answer is known.
A ring $R$ is called
 {\it Artin-Schelter {\em (}AS{\em )} Gorenstein}\label{asg-index} if 
$R$ has finite left and right injective dimension $d$ and, for some
shift $m$,
one has
$$\Ext^i(k,R) = \begin{cases} 0 & \text{if}\  i\not=d \\
k(m) & \text{if}\ i=d 
\end{cases}$$ 
Then \cite{YZ}  shows that
{\it  noetherian AS Gorenstein rings
of dimension $d+1$
satisfy the hypotheses of Theorem~\ref{YZ}.}  
 Moreover,  
$\omega^\circ$ is of the form  $\omega[d]$ where
$\omega$ is the image of an honest $R$-bimodule
 that is left and right free of rank one. 
It follows that, in the case of AS-Gorenstein rings,
Serre duality takes the form of natural isomorphisms
$$\Ext^i_{\qgr R}(-,  \omega)\cong
R^{d-i}\Gamma( -)^*.    $$
This is essentially the classical statement of Serre duality for
projective varieties \cite[Theorem~III.7.6]{H}.


\section{Artin-Schelter regular algebras}
\label{secasregular}

The area of mathematics described in this survey really began some
fifteen
years ago in the work \cite{AS} of Artin and Schelter, where they
attempted to
classify the appropriate analogues of a commutative polynomial ring in
three
variables. The final solution to this problem in \cite{ATV1,ATV2}  for
rings
generated in degree one and in \cite{Steph1,Steph2} for the general
case,
required the creation  of many of the ideas described in the earlier
sections
of this paper and so it  serves  as a convenient bridge between our
expositions  on noncommutative curves and on noncommutative surfaces.

The philosophy of Artin and Schelter is that a noncommutative projective
plane
should be represented by a ring with some of the most basic properties
of a
polynomial ring in three variables.  This point of view is
essentially adhoc but an a posteriori justification, which will be
discussed in
\S\ref{secncpq}, has been  provided by the work of Bondal and
Polishchuk in \cite{Bondal}.  In the approach of \cite{AS}, it is more
natural
to allow the generators of that polynomial ring to have varying weights,
but
this can be regarded as an advantage, since we will obtain (some)
noncommutative quadrics and more general weighted projective spaces by
the same
technique!

So 
 what should be the definition of our noncommutative analogue of a
polynomial ring in three variables? Obviously, it should be a connected
graded
ring, of homological 
dimension three and have finite Gelfand-Kirillov
dimension. By
themselves, these conditions are insufficient, since they do not exclude
some
particularly unpleasant rings like $k\{x,y\}/(xy)$. The solution was
to introduce   the Gorenstein condition mentioned in the last section.
 More precisely:

\begin{definition}\label{AS-regular}
A connected graded ring $A$ is called  {\em Artin-Schelter   (AS)
regular of dimension
$d$} if it satisfies the following conditions:
\begin{itemize}
\item[(1)]  $A$ has homological dimension $d$.
\item[(2)] $\GKdim(A)<\infty$. 
\item[(3)]  $A$ is  AS Gorenstein.
\end{itemize}
\end{definition}

The advantage of  the third hypothesis, 
for AS-regular rings of dimension $3$,
 is that the projective
resolution of $k$ is forced to be of the form 
\begin{equation}\label{gor-exact}
0\too A\too A^{n} \too A^{n} \too A\too k\too 0,
\end{equation}
(where  the shifts in degree have been ignored) for some $n$ and this 
quickly gives one strong information on the Hilbert series of $A$. Using
this,
 a complete classification was given for
AS regular rings of dimension three by Artin, Tate, Van den Bergh
and Stephenson:

\begin{theorem} \label{atv}
\cite{ATV1,ATV2, Steph1,Steph2}
The AS-regular rings $A$ of dimension $3$ can be classified. 
They are all noetherian domains
with the Hilbert series of a weighted polynomial ring 
$k[x,y,z]$---thus the $(x,y,z)$ can
be given  degrees $(a,b,c)$ other than $(1,1,1)$.
\end{theorem}

One of the more pleasant aspects of this classification is that
  the noetherian condition comes for free once one makes 
  the Gorenstein hypothesis.  
 Conversely, it could be that
the Gorenstein condition is automatic for  noetherian rings of finite
homological dimension \cite{StephZ}. 

In order to give an idea of the classification in Theorem~\ref{atv},
consider
the case when the variables all have weight $1$, so $A$ has Hilbert
series
$(1-t)^{-3}$.  (The proof in the general case is similar.)  Here, once
again,
one looks at $\qgr(A)$. The idea is to look at the simplest possible objects
in $\qgr(A)$, the (images of) {\it point modules}; 
modules $M=\bigoplus_{i\geq 0}M_i$ that are generated by $M_0$ and satisfy
$\dim_k M_i=1$, for all $i\geq 0$. If $A$ were commutative then such a module
would be isomorphic to a polynomial ring in one variable and hence be
isomorphic to the homogeneous coordinate ring 
of a point $p\in \Proj(A)$.
Hence the name. When $A$ is noncommutative it is easy to understand the point
modules combinatorially and this leads to the main step in the proof of 
Theorem~\ref{atv}: either $A\cong B(\PP^2, \Lscr, \sigma)$ for some $\Lscr$ and
$\sigma$ or  there exists a homomorphic image $B$ of $A$ isomorphic to
$B(E,\Lscr,\sigma)$ for some (possibly singular or
non-reduced) elliptic curve  $E$. Moreover, 
$B$ is defined by $3$ quadratic relations together 
with one cubic $\tilde{g}$
and so $A$ is defined by those $3$ quadratic relations.
 Notably, 
 $A$ is uniquely determined by $(E,\Lscr,\sigma)$!
Finally, one can show that image $g$ of $\tilde{g}$ in $A$ is 
a  non-zero divisor that is {\it normal} in the sense that 
it satisfies $gA=Ag$. Thus, $B = A/gA $
and this enables one to use 
one's detailed knowledge of the structure of $B$  
 to prove
the required properties for $A$.

It is clearest if we explain this procedure 
for an explicit example and we use 
\begin{equation}\label{silly1}
A\;=\; k\{x_1,x_2,x_3\}/(x_2x_1-\alpha x_1x_2,
x_3x_1-\beta x_1x_3,  x_3x_2-\gamma x_2x_3),
\end{equation}
where $\alpha, \beta,\gamma \in k^*$. The reader may recognize this ring as a
typical example of a ``quantum coordinate ring of affine $3$-space.''
Let $M=\bigoplus_{i\geq 0} m_ik$ be a point module for $A$, with 
module structure defined by $m_i x_j = \lambda_{ij}m_{i+1}$ for 
$\lambda_{ij}\in k$. Each relation  between the generators of
 $A$ must kill $m_0$ (or any $m_i$) and this leads to the following three
equations in the $\lambda_{ij}$:
$$\lambda_{02}\lambda_{11}- \alpha \lambda_{01}\lambda_{12}\;,\;
\lambda_{03}\lambda_{11}- \beta \lambda_{01}\lambda_{13}\;,\;
\lambda_{03}\lambda_{12}- \gamma \lambda_{02}\lambda_{13}.$$
These equations may be rewritten as
 $N\cdot \mathbf \lambda = 0$, where 
(ignoring the first subscript) 
$$N=N(\lambda) = \begin{pmatrix} \lambda_2 & -\alpha \lambda_1 & 0 \\
\lambda_3 & 0 & -\beta \lambda_1\\
0& \lambda_3 & -\gamma \lambda_2 \end{pmatrix}
\;{\rm and}\; 
\mathbf \lambda = \begin{pmatrix} \lambda_1 \\
\lambda_2 \\ \lambda _3\end{pmatrix}.$$
This  defines  a subvariety 
$\Lambda\subseteq \mathbb P^2\times \mathbb P^2
\equiv \mathbb P(A_1^*) \times \mathbb P(A_1^*).$
A little linear algebra shows that $\Lambda$ is 
the graph of an automorphism
$\sigma$ of a subvariety $E\subseteq \mathbb P^2$. 
Indeed, 
$$E=\{\mathbf \lambda \mid \det(N(\mathbf \lambda))=0\}
=\{(\lambda_1,\lambda_2,\lambda_3) \mid
(\beta-\alpha\gamma)\lambda_1\lambda_2\lambda_3 = 0\}.$$
Similarly, $\sigma$ is defined by $\sigma(\mathbf \lambda)
= \{\mathbf \mu : N(\mathbf \lambda ) \cdot \mathbf \mu=0\}.$ 
We leave it to the reader to check that $\mathbf \mu$ 
is uniquely determined and 
to determine the precise form of $\sigma$.
Since $E$ is determined by the point modules of $A$ it is often called the {\it
point scheme} $\Pscr(A)$\label{pointscheme-index} of $A$. 

Now form the twisted homogeneous coordinate ring
$B=B(E,\Oscr(1),\sigma)$. Then, $B_1=\HH^0(\Oscr(1))$ 
is naturally isomorphic to $A_1$ and 
it is an easy exercise to show that this
induces a homomorphism $\phi: A\to B$. The 
details can be found in 
\cite{ATV2} or \cite[Lemma 3.3]{SS}
but the idea is as follows: Given a
quadratic relation $f=\sum \alpha_{ij}x_ix_j$ 
we need to prove that $\phi(f)=0$.
By the definition of  multiplication in $B$ 
this amounts to proving that, if $\{\bar{x}_i : 1\leq i\leq 3\}$
denote the coordinate functions on $\mathbb P^2$,
then
\begin{equation}\label{silly}
\sum \alpha_{ij}\bar{x}_i(p)\bar{x}_j(p^\sigma) =
0\qquad{\rm  for\ all\ } p\in E.
\end{equation} 
Now, $\Lambda$ was defined as the common zeros of the linearizations
$\tilde{f} 
=\sum \alpha_{ij}x_{0i}x_{1j}$ of  all such  relations $f$.
 Thus, $\tilde{f}$ is zero on $\Lambda$ which, since 
$\Lambda=\{(p,p^\sigma) : p\in E\}$, implies that 
$\tilde{f}$ is zero on $\{(p,
p^\sigma) : p\in E\}$. This is equivalent to \eqref{silly}.
Finally, an application of the Riemann-Roch theorem 
can then be used to prove that $B$ 
is generated by $B_1$ and hence that $\phi$ 
is surjective. Thus, $B$ is a homomorphic image of $A$. 
The remarkable thing about Theorem~\ref{atv} is that 
this simple approach actually works for all 
AS regular algebras of dimension three.

Note that we really have two different types of rings in
\eqref{silly1}. 
If $\beta=\alpha\gamma$ then $E\cong \mathbb P^2$ 
and indeed $A\cong B$. 
Thus, $\qgr(A)\cong \mathbb P^2$,
 by Theorem~\ref{AZZ}. If 
$\beta\not=\alpha\gamma$ then 
$E=\{\lambda_1\lambda_2\lambda_3=0\}$ is the
product of three lines  and 
one can check that $B=A/x_1x_2x_3A$.
Thus Theorem~\ref{AZZ} now says that
$\coh(E)\cong \qgr(B)$ embeds into $\qgr(A)$. 
A natural interpretation is to regard $E$ as embedded
in the ``noncommutative projective plane'' 
 $\qgr(A)$. One can extend this analogy by 
observing that the (suitably defined)
complement $\qgr(A)\smallsetminus E
\cong \mod(R)$, where $R= A[g^{-1}]_0$ and so 
$\mod(R)$ is a ``noncommutative affine space.''
We reiterate that this phenomenon is typical 
of AS regular algebras of dimension three and
 illustrates the assertion from
the introduction  that 
all known noncommutative projective surfaces 
contain an embedded commutative curve. We will return to this 
point of view in \S\ref{noncom-planes}.

 The most interesting examples of AS regular 
algebras of dimension 
three are
the so-called $3$-dimensional Sklyanin algebras. 

\begin{example} \label{skl}
The ``generic'' case, when $E$ is a  smooth elliptic curve and 
$\sigma$ is given by translation
under the group law, corresponds to 
 the {\em Sklyanin algebra}
$$\Skl_3(a,b,c) =
k\{ x_0,x_1,x_2\}/(ax_ix_{i+1} + 
bx_{i+1}x_i + cx_{i+2}^2 : i = 1,2,3 \mod 3),$$
where 
$(a,b,c)\in \PP^2\smallsetminus F$, for a (known)
finite set $F$. \end{example}

The reader may like to repeat the earlier analysis and show that 
the associated curve for $\Skl_3=\Skl_3(abc)$ is the 
smooth elliptic curve $$E:\quad  abc(x_1^3+x_2^3+x_3^3) 
- (a^3+b^3+c^3)x_1x_2x_3=0.$$
Thus, $E$ embeds into $\qgr(\Skl_3)$.
An interesting challenge is to find a direct 
algebraic proof that $\Skl_3$ has Hilbert series $(1-t)^{-3}$. 
The only known proofs use the geometric 
techniques of this survey.

The results on three dimensional AS-regular rings have been extended to a
certain extent to higher dimension. Space does not permit us to to discuss this
in any detail so we will suffice by saying that, well before Artin-Schelter,
Sklyanin had constructed some putative four-dimensional algebras (see
\cite{Skly1,Skly2} and Example~\ref{skl4}, below) which were then generalized
by  Odesskii and Feigin  \cite{OF1,OF2}.  These are  now called higher
dimensional Sklyanin algebras.  That these algebras have the expected
properties was shown in a series of papers \cite{LS,SS,SmSz,TVdB}. 
However, in contrast to Theorem~\ref{atv}, there seem to exist a great many
other  four-dimensional AS-regular rings, and so there is no classification in
sight.

It is customary to divide the class of four-dimensional AS-regular algebras $R$
(generated in degree one) into subclasses according to the point scheme
$\Pscr(R)\subseteq \PP^3$ parametrizing the point modules.
  The case where $\Pscr(R)$ contains a quadric is 
best understood \cite{Vancliff1,VV}.  Many examples \cite{LSV,Skly1,Staf3},
including the $4$ dimensional Sklyanin algebra,  have been constructed where
$\Pscr(R)$  is the union of an elliptic curve and a finite set of points.
However, it is also possible for  $\Pscr(R)$  to consist of just a finite
set---indeed this is presumably the generic case. Very little is known about
the general properties of this class. For example, all
known rings $R$ with a finite point scheme are finitely generated as modules
 over their center 
\cite{VVW}, but it seems hard to believe that this could be true in
general. 

\begin{remark} Bondal has proved (unpublished) that any Poisson bracket
on $\PP^3$
  has at least a one-dimensional zero-locus. Thus  the existence of
  four dimensional AS-regular rings with a finite number of points
  indicates that intuition obtained from the quasi-classical case is
  not always reliable.
\end{remark}

In the later sections we will use factor rings of the $4$-dimensional
Sklyanin
algebra as examples of quadric surfaces and so  we end this section by
collecting
together the facts we need. For more details about this algebra and its
rich
module structure, we refer the reader to \cite{LS,SS,SmSz} or the survey
\cite{Sm}.

\begin{example}\label{skl4}
The {\em $4$-dimensional Sklyanin algebra $\Skl_4$} is the $k$-algebra
with  $4$ generators $x_0,\dots ,x_3$ and the $6$ relations
$$x_0x_i-x_ix_0=\alpha_i(x_{i+1}x_{i+2}+x_{i+2}x_{i+1}),\quad 
x_0x_i+x_ix_0=
x_{i+1}x_{i+2}-x_{i+2}x_{i+1},$$ where the subscripts are taken to be
$\{1,2,3\}$ mod $3$ and the $\alpha_i$ satisfy $\alpha_1\alpha_2\alpha_3
+ 
\alpha_1+\alpha_2+\alpha_3=0$ and $\{\alpha_i\}\cap\{0,\pm
1\}=\emptyset$.
The significant properties are  that $S=\Skl_4$ is an AS
regular algebra of dimension $4$ that 
has a two-dimensional space of  central  homogeneous
elements
$\Omega\subset S_2$. The factor ring 
$S/(\Omega)$ is isomorphic to a twisted homogeneous
coordinate
ring $B(E,\Lscr,\sigma)$, where $E$ is a
 smooth elliptic curve, 
$\sigma$ is an automorphism given by 
translation but $\Lscr$ is now a line
bundle of degree $4$. Since this data actually 
determines the Sklyanin algebra,
we sometimes  
 write $\Skl_4=\Skl_4(E,\Lscr,\sigma)$. 
In fact, $\Skl_4(E,\Lscr,\sigma)$ is 
even determined up to isomorphism
by $(E,\sigma)$, since $\Lscr$ just 
describes the presentation of the
algebra as a factor of a free algebra.  

Each $\Omega_\lambda\in \PP(\Omega)$ 
generates a prime ideal of $\Skl_4$ and we write 
$R_\lambda=R_{\lambda}(E,\Lscr,\sigma)$\label{Rlambda-index} for the factor
$\Skl_4/\Omega_{\lambda}\Skl_4$. 
\end{example}

\section{Commutative surfaces}
\label{seccomsur} \subsection{Overview} \label{seccomsurov}
After the classification of curves, the next classical subject in
algebraic
geometry is the classification of surfaces. Unfortunately this is a much
more
difficult subject and the information is far less complete. In this
section we
briefly survey this theory while in the remainder of the paper we will
discuss
the extent to which those results generalize to the noncommutative
universe.

For the rest of the paper,
 a ``surface'' will be a smooth proper connected scheme of
dimension two over $k$.  It is
well-known that such a surface is automatically projective.

Since the classification of surfaces $X$ is too hard to tackle
directly, one first classifies them up to birational
equivalence---this just means that one classifies the function fields
$k(X)$. Of course, $k(X)$ is just an extension of $k$ of transcendence
degree two. If $k(X)\cong k(x,y)$, we call $X$ a {\it rational
  surface}. A convenient interpretation of birational equivalence is
given by the theory of minimal models; if $X$ is a surface then $X$ is
{\it minimal}\label{minimal-index}
 if every surjective map $X\r Y$ to a birational surface
$Y$ is an isomorphism.  Every surface is birationally equivalent to a
minimal one and this minimal model is unique except for ruled and
rational surfaces. For rational surfaces the minimal models are
respectively $\PP^2$ and $\PP^1$-bundles over $\PP^1$ 
(these are also known as  Hirzebruch surfaces) 
 \cite[Thm
V.10]{Beauville}.  For more details, we refer the reader to
\cite{BPV,Beauville,H}.

By contrast, the classification  of  surfaces
within a birational equivalence class has a much neater answer. A
classical result by Zariski states that if $X$ and $X'$ are surfaces
which are birationally equivalent then it is possible to go from $X$
to $X'$ by repeatedly  applying blowing-ups and blowing-downs.

As we will want to mimic the procedure later, 
we first  describe blowing up/down in more detail (see, also,
 \cite[\S V.3]{H}). Let
 $p\in X$ be a point on the surface $X$,
with corresponding maximal ideal
$\mathfrak{m}$. Then the blow-up $\widetilde{X}$ of  $x$ in $X$  is defined as $\Proj
\Dscr$
where 
\begin{equation}
\label{standardreesring}
\Dscr\;=\;\Oscr_X\oplus \mathfrak{m}\oplus \mathfrak{m}^2\oplus\cdots,
\end{equation}
thought of as a sheaf of $\Oscr_X$-algebras on $X$ (for the
construction of $\Proj$ of a sheaf of graded algebras see \cite{H}). 
Since blowing up is
an essentially local transformation, it is enough to analyse this in the
case
when  $X=\Spec k[x,y]$ and  $p$ is the origin.  After some
elementary computations one obtains
\[
\widetilde{X}=\Spec k[y/x,x]\cup \Spec k[x/y,y]
\]
It follows that, away from $p$, the map $\widetilde{X}\r X$ is an
isomorphism.
However the inverse image $L$ of $p$ is isomorphic to $\PP^1$ and
has self intersection $-1$ \cite[Proposition~V.3.1]{H}. The curve $L$
 is called the \emph{exceptional curve}.\label{excep-curve} More
generally,
if $Y$ is an arbitrary surface     then a curve $L\subset Y$ is called
exceptional
if $L\cong \PP^1$ and if $L$ has self intersection~$-1$.

The inverse to blowing up is blowing down for which, given  an exceptional
curve $L$ on a surface $Y$, one constructs a  map (determined up to a
unique isomorphism) $\alpha:Y\r X$ to
some surface $X$ such that $\alpha$ contracts $L$ to a point and  is an
isomorphism outside $L$. One way of constructing $\alpha$ is outlined below, 
and this technique also provides some  birational transformations of a more
global nature than just blowing up or down.

Thus, 
 let $\Lscr$ be a line bundle on a surface  $X$ and fix nonzero  
global sections $x_0,\ldots,x_n\in \text{H}^0(X,\Lscr)$. If $\Lscr$ is 
generated  by the  $x_i$, then they define 
a map $x:X\r \PP^n$  such that $\Lscr=x^\ast \Oscr(1)$ and such that the $x_i$
 are the pullbacks of the coordinate functions \cite[Theorem
II.7.1]{H}.  Let $X'$ denote the
 image of $X$. In general there is not much to
be said about $X'$ although, in some very nice cases, $x$ will just contract a
finite number of exceptional curves. In particular,  if $L$ is an
exceptional curve on $X$, then one can pick $\Lscr$ such that $x$ contracts
precisely $L$ and nothing else \cite[proof of Thm V.5.7]{H}. This is the 
classical construction of blowing down.

Suppose now that $\Lscr$ is not generated by $x_0,\ldots,x_n$ and
let $\Fscr$ be the subsheaf of $\Lscr$ generated by those sections.
In this case, it is convenient to  replace   $\Lscr$ by 
$\Fscr^{\ast\ast}$ so that we may assume that $\Lscr/\Fscr$ has 
finite length.     In other words,
if
\[
U=\{p\in X\mid \Lscr\text{ is generated by $x_0,\ldots,x_n$ in $p$}
\}
\]
then $X\setminus U$ is a finite set of points and $x$ defines a map
of algebraic varieties $x:U\r \PP^n$. Classically, one calls this a
{\it rational map} $x:X\r \PP^n$.  We define the
 image $X'$ of $x$ to be the Zariski-closure of $x(U)$ in $\PP^n$.
 
In the following three subsections we discuss some classical birational
transformations for which we will be able to construct noncommutative
analogues
in \S\ref{sec12}. 

\subsection{Cubic surfaces \cite[Proposition 4.8]{H}}
\label{seccub}
Let $X=\PP^2$ and let $p_1,\ldots,p_6\in \PP^2$ be six points in
general position. By this we mean that no two points are equal, no
three lie on a line and not all six lie on a conic. Then the space of
sections of $\Oscr(3)$ (cubic forms) vanishing in $p_1,\ldots,p_6$ is
four dimensional and these sections generate $\Oscr(3)$ on
$\PP^2\smallsetminus\{p_1,\ldots,p_6\}$. 
As explained above we obtain a corresponding
rational map $x:\PP^2\r \PP^3$. In this case it turns out that the
image of $x$ is a cubic surface in $\PP^3$ and it is obtained from
$\PP^2$ by blowing up the six points $p_1,\ldots,p_6$. It is a
classical fact \cite[\S 24]{Manin2} that any cubic surface in $\PP^3$
can be obtained in this way.

\subsection{The Cremona transformation \cite[Example V.4.2.3]{H}}
\label{seccrem}
Fix three points $p,q,r \in \PP^2$ that are not
collinear. The space of sections of $\Oscr(2)$ (quadratic forms)
vanishing in those three points is three dimensional and so we now
obtain a rational map $x:\PP^2\r \PP^2$. In this case $x$ is 
birational and is classically called a  {\it Cremona transformation}.

If we choose homogeneous coordinates $x,y,z$ and we put $p,q,r$
in the standard positions $(1,0,0)$, $(0,1,0)$, $(0,0,1)$ then this
Cremona transformation  is given by $(x,y,z)\mapsto (yz,zx,xy)$. From
this explicit description it follows  easily  that the Cremona
transformation is the composition the following birational
transformations.
\begin{itemize}
\item[(a)] First   blow up the points $p,q,r$ to obtain a
  surface $\widetilde{X}$.
\item[(b)] It follows from general theory that the three lines
  connecting the points $p,q,r$ become exceptional and disjoint
  in $\widetilde{X}$. Now blow down these exceptional lines.
\end{itemize}
\subsection{Quadrics \cite[Exercise V.4.1]{H}} 
\label{secquad}
Now we choose two distinct points $p,q\in \PP^2$ and we look at the four
dimensional space of sections of $\Oscr(2)$ which vanish in
$p,q$. This gives us a rational map $x:\PP^2\r \PP^3$ whose image is
a quadric surface. In this case $x$ is obtained from blowing up
$p,q$ and then blowing down the connecting line.

The inverse to $x$ (as a birational map) can be 
constructed as follows.   Every quadric is
isomorphic to $\PP^1\times \PP^1$. Fix some point $r\in
\PP^1\times \PP^1$ and consider the three dimensional space of
sections of $\Oscr(1,1)$ vanishing in $r$. This gives us a
rational map $y:\PP^1\times \PP^1\r \PP^2$ which is the inverse to
$x$.

\subsection{Some reformulations}
\label{secsomref}
Unfortunately the constructions exhibited in  \S\ref{seccomsurov} do
not generalize well to the noncommutative case and so 
 we will give a more suitable reformulation. 

Let $X$ be a surface and let $\Lscr$ be a line bundle on $X$ generated by
global sections. Then $\Lscr$ defines a map $x:X\r\PP^n$ and in
\S\ref{seccomsurov} we have considered $X'=x(X)$.  With an eye towards
noncommutative generalizations, it is better to consider  
$\widetilde{X}=\Proj A$,  where $A$ is the ring $\bigoplus_n
H^0(X,\Lscr^{\otimes n})$.  To understand the relation between $X'$ and
$\widetilde{X}$ we introduce the Stein factorization 
$x:X\xrightarrow{\alpha}
Y\xrightarrow{\beta} \PP^n$ of $x$. This factorization has the defining
property that $\beta$ is finite and $\alpha_\ast(\Oscr_X)=\Oscr_Y$. Put
$\Lscr_Y=\beta^\ast \Oscr(1)$.  This yields 
$$ \begin{array}{rl}
H^0(X,\Lscr^{\otimes n})&= H^0(X,\alpha^\ast(\Lscr_Y)^{\otimes
n})
=H^0(Y,\alpha_\ast\alpha^\ast(\Lscr_Y)^{\otimes
n})\\ 
\noalign{\vskip 2pt}
& =H^0(Y,\alpha_\ast \Oscr_X \otimes_{\Oscr_Y} \Lscr_Y^{\otimes
n})
=H^0(Y,\Lscr^{\otimes n}_Y).
\end{array}
$$
  Since $\beta$ is finite, $\Lscr_Y$ is ample
and so, by Serre's Theorem~\ref{serre}, we obtain $\widetilde{X}=Y$.  Hence
$\widetilde{X}$ is a finite cover of $X'$.

If $\Lscr$ is not generated by its global sections then,  as in
\S\ref{seccomsurov}, we let $\Fscr$  be the subsheaf of $\Lscr$ spanned
by the
global sections of $\Lscr$ and we define $\widetilde{X}=\Proj A$ for
$A=\bigoplus H^0(X,\Fscr^n)$. Here $\Fscr^n$ is, by definition, the
image of
$\Fscr^{\otimes n}$ in $\Lscr^n$. As in the previous paragraph, one can
show
that $\widetilde{X}$ is a finite cover of $X'$.

In the standard examples that were considered in
\S\S\ref{seccub}--\ref{secquad},   $\widetilde{X}=X'$.

\section{Noncommutative surfaces}\label{secnoncom-surface}
We start this section with the embarrassing admission that we do not
really know what the correct definition should be for a
noncommutative surface. Keeping with the philosophy of 
\S\S\ref{sec4}--\ref{secy1} it would be tempting to accept one of 
the following
definitions, the first of which is more algebraic while the second,
inspired by
\cite{Bondal1}, is more geometric.
\begin{definition}\label{attemptI}(Attempt I) A noncommutative 
projective (normal but possibly singular)
  surface $\Cscr$ is a category of the form $\QGr(A)$ where $A$ is a
  noetherian connected
  graded domain  of Gelfand-Kirillov dimension
  three such that $A$ is a maximal order.  Moreover, $\Cscr$ is smooth
if 
 $\QGr(A)$ has homological
  dimension $2$.
\end{definition}

\begin{definition} (Attempt II) 
\label{attemptII} A smooth  noncommutative surface is a
  locally noetherian Grothendieck category of homological dimension
  two such that the full subcategory of noetherian objects is
  $\Ext$-finite and saturated.
\end{definition}

We should make some comments about Definition \ref{attemptI}.
Probably one should work with prime rings, but 
this introduces technical complications unrelated to the 
underlying intuition, as was true in \S\ref{sec5} 
 vis-a-vis \S\ref{sec4}. The assumption that 
  $A$ is a maximal order is natural, given the discussion in
\S\ref{secmorecomments}.   
Finally, this definition tacitly assumes the conjecture (see
\cite[p.233]{Staf5}) that $2<\GKdim(A)<3$ is impossible for a connected
graded domain $A$. By \cite[Theorem 0.1]{Staf5} it is known that 
$2<\GKdim(A)<2\frac{1}{5}$ is impossible.

Unfortunately, at this point it seems to be completely intractable to
obtain any
substantial classification results starting from these definitions, so
the
precise definition is moot.  Instead, the study on noncommutative
surfaces has
been more involved with an attempt to extend some basic parts of the 
commutative theory. This is what we will discuss in the forthcoming
sections. In particular, as we will describe, there are good arguments 
that
all noncommutative $\PP^2$'s and quadrics have been 
described and that, at least for
important special cases,  we have the right analogue of blowing up and
down.

\subsection{Birational classification}\label{birat-class}
Mimicking the commutative theory we should first obtain a birational
classification for noncommutative surfaces. However there are
immediately some foundational questions here.

If we accept  Definition \ref{attemptI}, then  
the corresponding function field 
is $Q(A)_0$, in the notation of \S\ref{sec4}. 
For Definition \ref{attemptII}   the function field should presumably 
be defined as the
endomorphism ring of the direct sum of the maximal injectives, but it
is not clear that this notion is well behaved.
Alternatively, we can concentrate
  initially  on the classification of division algebras $D$
  of transcendence degree two. Here, we use the definition of
transcendence
  degree from \cite{GK} (essentially, this means that $D$ is the Goldie
  quotient ring of a finitely generated, 
  (ungraded) domain $R$ of Gelfand-Kirillov dimension $2$). 
 Hopefully this is equal to the more sophisticated
  definition from \cite{Zhang2}, for which more favorable properties
  have been proved. For other types of transcendence degree of a more
  homological nature see \cite{Resco,Schofield}. 

If we try to write down a list of known division algebras of
transcendence degree two then  the list is surprisingly short:
\begin{itemize}
\item[(a)] {\it Division algebras which are finite over a central
  commutative subfield of transcendence degree two.}
\item[(b)] {\it Division rings of fractions  of 
skew polynomial extensions of $k(X)$, for curves $X$.} 
\cite[p.439]{Cohn}
\item[(c)] {\it The function fields $Q(\Skl_3)_0$
of the three dimensional Sklyanin
algebras $\Skl_3$.}
\end{itemize}
Since no other examples are known, Mike Artin makes 
the provocative conjecture in \cite{Ar2} that these are in fact all the
division
algebras of transcendence degree two. The first author of this survey
feels that this conjecture is false, whereas the
second author, being more diplomatic, is neutral!

Anyway, true or not, this conjecture immediately gives rise to some
interesting algebraic problems:

\begin{examples} 
\label{moredetailsfollow}
Let $E$ be the function field of the ``affine quantum plane''. That
is, $E$ is generated as a division algebra by two variables $u,v$
satisfying the relation $vu=quv$ for some $q\in k^\ast$ which is not a
root of unity. Let $\tau$ be the automorphism of $E$ given by
$u\mapsto u^{-1}$, $v\mapsto v^{-1}$ and put $D=E^\tau$.  Then $D$ is
a division algebra of transcendence degree two and the question is
where it fits in the above list.  It turns out that $D$ is also the
function field of a two-dimensional affine quantum plane and this was first
proved using the techniques outlined in this survey. The details will be given
in \S\ref{secmoredetails}. 
\end{examples}

\section{Noncommutative projective planes and quadrics}
\label{secncpq}
Starting with this section we will make our first steps into the
realm of noncommutative surfaces. In order to motivate the
definition of the two most basic examples, namely the noncommutative
analogues of $\PP^2$ and $\PP^1\times \PP^1$ we discuss an extension
of Theorem \ref{AZ}. This material was to some extent inspired by
\cite{Bondal}. 

\subsection{Extension of the Artin-Zhang Theorem}
\label{sec8}
Theorem \ref{AZ} depends on the availability of a suitable
autoequivalence $s$ which sometimes makes it hard to apply. 
 Inspection of the condition for ampleness of $s$ reveals
that it only depends on the objects $\Oscr(n)=s^n \Oscr$ and not on $s$
itself. On the other hand, the construction of the ring $B$ depends on
$s$. If we are willing to replace the notion of a $\ZZ$-graded
$k$-algebra by the   weaker notion of a $\ZZ$-algebra, then we
can get
rid of $s$ altogether. In general a $\ZZ$-algebra is just a $k$-linear category
whose objects are indexed by $\ZZ$ and a module over a $\ZZ$-algebra is
a $k$-linear functor to the category of $k$-vector spaces. The idea of
regarding
additive categories as generalizations of rings is  standard. See for
example \cite{Mi}.

Let us now spell out the abstract definition of a $\ZZ$-algebra in
concrete terms. 
 In the sequel a \emph{$\ZZ$-algebra}\label{Z-algebra-index}
  \cite{Bondal} will be a $k$-algebra
$A$
(usually without unit) together with a decomposition $A=\bigoplus_{ij\in
\ZZ}
A_{ij}$ such that the multiplication has the property that
$A_{ij}A_{jk}\subset
A_{ik}$ and $A_{ij}A_{kl}=0$ if $j\neq k$.  It is clear that the
$A_{ii}$ are
subalgebras of $A$. We require that these have units $e_{i}$ such that
if $a\in
A_{ij}$ then $e_{i}a=a=ae_{j}$. Following \cite{Bondal} our point of
view will
be that $\ZZ$-algebras are generalizations of $\ZZ$-graded rings in the
following sense: if $B=\bigoplus_{i\in \ZZ}B_i$ is a $\ZZ$-graded
$k$-algebra
then putting $A_{ij}=B_{j-i}$ defines a $\ZZ$-algebra. Most properties
of
graded rings can be stated and proved in the more general context of 
$\ZZ$-algebras.

Given a $\ZZ$-algebra $A$, then  a right $A$-module will be an ordinary
right
$A$-module $M$ together with a decomposition $M=\bigoplus M_i$ such that
$M_i
A_{ij}\subset M_j$, the $e_i$ act as units on $M_i$ and $M_i A_{jk}=0$ if
$i\neq
j$. We denote the category of right $A$ modules by $\Gr(A)$.
 It is clear that $\Gr(A)$ is a
Grothendieck
category and that, if $A$ is obtained from a $\ZZ$-graded ring $B$, then
$\Gr(A)\simeq\Gr(B)$. We will say that $A$ is noetherian if $\Gr(A)$ is
a locally
noetherian category.  We  borrow other notions from the theory of graded
rings.
For example, $M\in \Gr(A)$ is called \emph{right bounded} if $M_i=0$ for
$i\gg
0$. We say that $M$ is \emph{torsion} if it is a direct limit of right
bounded
objects.  We denote the corresponding category by $\Tors(A)$ and set
$\QGr(A)=\Gr(A)/\Tors(A)$.

Now consider  a locally noetherian Grothendieck category $\Cscr$ and let
$({{\Oscr}}(n))_{n\in\ZZ}$ be a series of noetherian objects in $\Cscr$.
Associated to the ${{\Oscr}}(n)$ there is a  $\ZZ$-algebra $A$ given by
$A_{ij}=\Hom({{\Oscr}}(-j),{{\Oscr}}(-i))$ for $j\ge i$, $A_{ij}=0$
otherwise,  and multiplication defined by composition.
We say that $({{\Oscr}}(n))_n$ is \emph{ample} if the conditions of 
Definition~\ref{ampleness} hold. 

We have the following obvious analogue of Theorem~\ref{AZ}.

\begin{theorems} \label{AZ1} Keep the above notation and 
assume that $\Cscr$ has finite
  dimensional $\Hom$'s between noetherian objects. If $({{\Oscr}}(n))_n$
is
  ample then  $ \Cscr\simeq\QGr(A) $. In
  addition $A$ is noetherian. 
\end{theorems}

\subsection{Noncommutative planes}\label{noncom-planes}
Following Artin and Schelter we propose the following definition of a
noncommutative $\PP^2$. It is somewhat more convenient to work
with $\QGr$ rather than $\qgr$ but, by \S\ref{notation}, the difference is not
significant.  Although $\ZZ$-algebras do not appear in its
definition, they will appear in its justification. 

\begin{definitions}\label{noncom-p2}
  A {\em noncommutative $\mathbb P^2$} is a Grothendieck 
  category of the form $\QGr(A)$,  where $A$ is an 
AS regular algebra of dimension three, with Hilbert series $(1-t)^{-3}$.
 \end{definitions}

Let $A$ be such  an 
AS regular algebra.
{}From the discussion in \S\ref{secasregular}, either 
$A\cong B(\PP^2,\Oscr(1),\sigma)$  for some $\sigma$---in which case 
Theorem \ref{AZZ} implies that $\QGr(A)\cong \PP^2$---or
$A$ has a normal element $g$ such that 
 $A/gA$ is isomorphic to $B(E,\Lscr,\sigma)$ for some 
elliptic curve $E$. In  the latter case we call $\QGr(A)$ a
 {\it nontrivial noncommutative $\PP^2$}\label{noncom-p22}.

We emphasize that this shows that $A$ is {\it not} determined by 
$\QGr(A)$. However, when one works with 
$\ZZ$-algebras, as we do below, this issue goes away since
 Corollary~\ref{noncanonical} shows that the 
$\ZZ$-algebra associated to $\QGr(A)$ 
is determined, up to isomorphism, by that category.

One may be tempted to extend Definition \ref{noncom-p2} to say that a 
noncommutative
weighted $\mathbb P^2$  is the category $\QGr(A)$, for an arbitrary
AS regular algebra $A$ of dimension three. However, as we show in
the next section, this is too restrictive.  But, if that is the case, is
the 
definition of a  noncommutative $\mathbb P^2$ also too restrictive?  As
an
indication that the definition is appropriate, we show that any possible
deformation (we use the term only in an informal sense that will become
clearer
as we proceed) of $\mathbb P^2$ will fit into the above classification,
and so
this suggested definition of a noncommutative $\PP^2$ is appropriate.

The approach we use for $\PP^2$ was found
by Bondal and Polishchuk in \cite{Bondal} and  we give a slightly
personalized
interpretation of the material in that paper. 
Let $\Cscr$ be our supposed noncommutative projective plane. 
 The point of view  of  \cite{Bondal} is that a
good definition for a noncommutative $\PP^2$ should at least
encompass deformations of $\Qch(\PP^2)$ as a category. 
 One way to control deformations of
 commutative schemes is through the structure of their
exceptional objects. Recall that an object $E$ in an abelian category
is \emph{exceptional}\label{excep-ob-index} if it satisfies
\begin{equation}
\label{exceptionalobject}
\Ext^i(E,E)=
\begin{cases}
k&\text{if $i=0$}\\
0&\text{otherwise}
\end{cases}
\end{equation}
Standard deformation theory yields that exceptional 
objects lift uniquely to a deformation.

On $\PP^2$ the line bundles $\Oscr_{\PP^2}(n)$ are clearly exceptional
objects
and so it is reasonable to impose the existence of analogous objects
${{\Oscr}}(n)$ in $\Cscr$. Now the ampleness of a sequence is also
stable under
deformation \cite{VdB26}, so we will assume  that the ${{\Oscr}}(n)$ in
fact
form an ample sequence in $\Cscr$. By Theorem \ref{AZ1}, $\Cscr$ will be
equivalent to $\QGr(A)$  where $A$ is the $\ZZ$-graded algebra defined by 
  $A=\bigoplus_{j\geq i}A_{ij}$ with
$A_{ij}=\Hom({{\Oscr}}(-j),{{\Oscr}}(-i))$. Note in passing that we
do not
have a natural $\NN$-graded algebra associated to our data. The problem
is that
if we, say, took  $A'=\bigoplus_{i\geq
0}\Hom({{\Oscr}}(0),{{\Oscr}}(-i))$,
then we have no convenient shift functor and so it is unclear how one
defines
the multiplication on $A'$.

To conclude we  have to describe $A$ more precisely. 
The main idea is that the dimension of
$\Ext^t_\Cscr({{\Oscr}}(-j),{{\Oscr}}(-i))$ may
decrease under deformation but  the Euler form
\[
\langle {{\Oscr}}(-j) ,{{\Oscr}}(-i)\rangle=
\sum_t (-1)^t\dim\Ext^t({{\Oscr}}(-j),{{\Oscr}}(-i))
\]
should remain constant. This allows us to compute $\dim A_{ij}$ by
using the fact that
$\Ext^t_{\Oscr_{\PP^2}}(\Oscr_{\PP^2}(-j),\Oscr_{\PP^2}(-i))=0$ for
$t>0$ and
$j\ge i$. Putting this together with a little linear algebra shows that 
$A$ must satisfy the   properties of the following definition:
 
\begin{definitions} 
\label{noncommutativeP2}
A  (noncommutative) {\em Bondal-Polishchuk $\PP^2$} is an abelian
category
  of the form $\QGr(A)$ where $A$ is a $\ZZ$-algebra
  with the 
following properties:
\begin{itemize}
\item[(P1)] $A$ is generated by the $A_{i,i+1}$.
\item[(P2)] The following holds for the dimensions of $A_{ij}$:
\[
\dim A_{i,i+n}=
\begin{cases}
(n+1)(n+2)/2 &\text{if $n\ge 0$}\\
0&\text{otherwise}
\end{cases}
\]
\item[(P3)] Put $V_i=A_{i,i+1}$. By (P2), the multiplication map
  $V_i\otimes V_{i+1}\r A_{i,i+2}$ has a $3$-dimensional
  kernel  $R_i$. We require that the $R_i$
  generate the relations of $A$.
\item[(P4)] Put $W_i=V_i\otimes R_{i+1}\cap R_i\otimes V_{i+2}$ inside
 $ V_i\otimes V_{i+1}\otimes V_{i+2}$. Dimension
  counting using (P2) and (P3) reveals that $\dim W_i=1$. 
  If  $W_i=w_ik$, we require that $w_i$ is a
  non-degenerate tensor, 
   both as an element of $V_i\otimes R_{i+1}$ and as an
  element  of $R_{i}\otimes V_{i+2}$.
\end{itemize}
\end{definitions}
It is easily seen that if $A$ is a three dimensional 
AS regular algebra with Hilbert series $(1-t)^{-3}$
and associated
$\ZZ$-algebra $A'$, then $A'$ satisfies the above axioms.

In  \cite{Bondal} the authors classify the  $\mathbb Z$-algebras 
satisfying
(P1)--(P4) in terms of a triple $(Y,\Lscr_1,\Lscr_2)$, where $Y$ is
either
$\PP^2$ or an elliptic curve therein and the $\Lscr_i$ are ample line
bundles
on $Y$ satisfying certain technical  conditions. These imply that
$\Lscr_2=\Lscr_1^\sigma$ for an automorphism $\sigma$ on $Y$. Thus
 it is not
surprising that they further obtain:

\begin{theorems}
\label{moregeneral} Assume that $\charact k\not= 3$.
 Then the  noncommutative {\it Bondal-Polishchuk $\PP^2$}'s
are exactly the noncommutative $\PP^2$'s. More precisely,
every $\ZZ$-algebra satisfying (P1)-(P4) is isomorphic
   to a $\ZZ$-algebra  obtained from a connected graded 
   AS algebra of dimension three with Hilbert
   series $(1-t)^{-3}$.
\end{theorems}

This theorem  provides another  justification for our definition of a
noncommutative $\PP^2$. In \cite{Bondal} it was only stated in
characteristic zero.

If $A$ is a $\ZZ$-algebra and $m\in \ZZ$ the we may
define the {\it shifted $\ZZ$-algebra} $A(m)$ by $A(m)_{ij}=A(m)_{i+m,j+m}$.
Clearly $\Gr(A(m))\cong \Gr(A)$ and  $\QGr(A(m))\cong\QGr(A)$ but
  $A$ comes from a graded ring if and
only if $A\cong A(1)$.
If $\Cscr$ is a  noncommutative Bondal-Polishchuk $\PP^2$ then it is easy to see that
  the  sequence of objects $(\Oscr(n))_n$ is determined up to shift by the
categorical  properties of $\Cscr$ and hence that
 the associated $\ZZ$-algebra is also
  determined up to shifting (see \cite{VdB26} for the details).
  Combined with Theorem~\ref{moregeneral} these observations prove:
  
\begin{corollarys}\label{noncanonical}
Assume that 
$\operatorname{char} k\neq 3$. If  $\Cscr$ is a noncommutative $\PP^2$ then the associated 
$\ZZ$-algebra $A$ is unique. For example, if $\Cscr=\Qch(\PP^2)$,
then $A$ is just the $\ZZ$-algebra associated with $k[x,y,z]$.
\end{corollarys}

\subsection{Noncommutative quadrics}
\label{sec10}

Even if one restricts to algebras generated in degree one, there are two
classes of Artin-Schelter regular algebras of dimension three; the ones
already
mentioned with Hilbert series $(1-t)^{-3}$  and
a second class $\Qscr$  of algebras that require
 just two generators and have the Hilbert
series $(1-t)^{-2}(1-t^2)^{-1}$ of a weighted polynomial ring
$k[u_1,u_2,u_3]$
with weights $(1,1,2)$.  The archetypal example of this second class is
the
enveloping algebra $U(\mathfrak g)$ of the three-dimensional Heisenberg
Lie
algebra $\mathfrak g$; thus $\mathfrak g$ has basis $x,y,z$ with
relations
$[x\,y]=z,$ $[x\,z]=0=[y\,z]$.  Clearly $U(\mathfrak g)$ is generated by
$x$
and $y$, and if we assign $x,y,z$ weights $1,1,2$, then the
Poincar\'e-Birkoff-Witt Theorem ensures that $U(\mathfrak g)$ is a
graded ring
of the required Hilbert series.  Given the results from the last
section,
one might propose that a noncommutative quadric surface be defined as
$\QGr(A)$
for any $A\in \Qscr$.  (These surfaces are smooth---see
Remark~\ref{serre2}.)

Unfortunately,  this is inadequate since 
there appear to 
be far   more noncommutative quadrics than AS regular
rings. For example, the Hilbert series of the factor
 ring $R_\lambda=R_{\lambda}(E,\Lscr,\sigma)$ of the $4$
dimensional Sklyanin algebra, as described in Example~\ref{skl4}
is the same as that of the homogeneous coordinate ring of a quadric in
$\PP^3$. Therefore we would like to consider $R_\lambda$ as a
noncommutative quadric (see below and \S\ref{sec11} for a more  formal
justification). But $R_\lambda$ depends upon $3$ paramenters
(one each for $E$, $\Lscr$, $\sigma$ and $\lambda$ minus one for the
automorphism group of $E$)
whereas a similar count reveals that the   AS regular
rings depend on just  $2$.

In order to get an idea about how many noncommutative quadrics one
should
expect, we follow the Bondal-Polishchuk approach from the last section,
but
with $\PP^2$ replaced by $\PP^1\times \PP^1$. 
For this variety some of the exceptional line bundles are:
\begin{equation}\label{exceplb}
\Oscr(n)=
\begin{cases}
\Oscr(k,k)&\text{if $n=2k$}\\
\Oscr(k,k+1)&\text{if $n=2k+1$}
\end{cases}
\end{equation}
where $\Oscr(a,b)=\Oscr_{\PP^1}(a)\boxtimes \Oscr_{\PP^1}(b)$. If one
now
follows
the approach outlined in the last section one arrives at:

\begin{definitions}\label{noncom-quadrics}
\cite{VdB26} {\em A noncommutative quadric} is a Grothendieck
  category of the form $\QGr(A)$ where $A$ is a $\ZZ$-algebra satisfying
the following properties.
\begin{enumerate}
\item[(Q1)] $A$ is generated by the $A_{i,i+1}$.
\item[(Q2)] The following holds for the dimensions of $A_{ij}$
\[
\dim A_{i,i+n}=
\begin{cases} 
0&\text{if $n<0$}\\
(k+1)^2&\text{if $n=2k$ and $n\ge 0$}\\
(k+1)(k+2)&\text{if $n=2k+1$ and $n\ge 0$}
\end{cases}
\]
\item[(Q3)] Put $V_i=A_{i,i+1}$. By (Q2), $A_{i,i+2}=V_i\otimes
  V_{i+1}$ and the multiplication map $V_i\otimes V_{i+1}\otimes
  V_{i+2}\r A_{i,i+3}$ has two-dimensional kernel
   $R_i$. We require that the $R_i$ generate the relations of
  $A$.
\item[(Q4)] Put $W_i=V_i\otimes R_{i+1}\cap R_i\otimes V_{i+3}$ inside
  $V_i\otimes V_{i+1}\otimes V_{i+2}\otimes V_{i+3}$. Dimension
  counting using (Q2) and (Q3) reveals that $\dim W_i=1$. If $W_i=w_ik$, 
  we require that $w_i$ is a
  non-degenerate tensor, both as an element of $V_i\otimes R_{i+1}$ and 
 as an element of  $R_{i}\otimes V_{i+3}$.
\end{enumerate}
\end{definitions}

As yet there is no analogue of Theorem \ref{moregeneral}, but the
$\ZZ$-algebras satisfying (Q1)-(Q4) are surely more general that the
AS-regular algebras $A\in\Qscr$, simply because
 they depend on one more parameter. 
In more detail, the classification in \cite{VdB26} 
of $\ZZ$-algebras satisfying (Q1)-(Q4) 
is in terms of quadruples
$(Y,\Lscr_0,\Lscr_1,\Lscr_2)$ where either:
\begin{enumerate}
\item $(Y,\Lscr_0,\Lscr_1,\Lscr_2)=
(\PP^1\times \PP^1,\Oscr(1,0),\Oscr(0,1),\Oscr(1,0))$;
  or else
\item \begin{enumerate}
\item $Y$ is a
curve which is embedded as a divisor of degree $(2,2)$ in $\PP^1\times
\PP^1$ by the global
sections of both $(\Lscr_0,\Lscr_1)$ and $(\Lscr_1, \Lscr_2)$.
\item  $\Lscr_0\not\cong \Lscr_2$.
\item  $\deg \Lscr_0\mid
C=\deg\Lscr_2\mid C$ for every irreducible component $C$ of $Y$. 
\end{enumerate}
\end{enumerate}
The construction of the quadruple associated to a $\ZZ$-algebra $A$ 
is similar to the Bondal-Polishchuk approach in \cite{Bondal}: we
consider $R_0$ as the equations for a closed subscheme $Y$ of
$\PP(V_0^\ast)\times \PP(V_1^\ast)\times \PP(V_2^\ast)$ and we put
$\Lscr_i=\pr^\ast_i(\Oscr(1))$ where $\pr_i$ denotes the projection
onto $\PP(V_i^\ast)$. Conversely, to construct a $\ZZ$-algebra from a
triple $(Y,\Lscr_0,\Lscr_1,\Lscr_2)$, one may proceed as follows. First
construct a sequence of line bundles 
(called  an elliptic helix in \cite{Bondal})
$(\Lscr_i)_i$ defined by the property
\begin{equation}
\label{sequence}
\Lscr_i\otimes \Lscr_{i+1}^{-1}\otimes 
\Lscr_{i+2}^{-1}\otimes \Lscr_{i+3}=\Oscr_Y.
\end{equation}
We put $V_i=\Gamma(Y,\Lscr_i)$ and  let $R_i$ be the kernel of the
canonical map $$\Gamma(Y,\Lscr_i)\otimes \Gamma(Y,\Lscr_{i+1})\otimes
\Gamma(Y,\Lscr_{i+2})\r \Gamma(Y,\Lscr_i\otimes \Lscr_{i+1}\otimes
\Lscr_{i+2}).$$
 Finally, $A$ is defined to be the $\ZZ$-algebra freely generated
by the $V_i$, subject to the relations $R_i$. 

An interesting open question is to determine when $\QGr(A)\cong
\QGr(A')$ for
$\ZZ$-algebras satisfying (Q1)-(Q4). This is even unanswered for
the 
factor rings $R_{\lambda}=\Skl_4/(\lambda_1\Omega_1+\lambda_2\Omega_2)$ 
of the
Sklyanin algebra $\Skl_4$ described in Example~\ref{skl4}. (As we will see in
the next section, generically these rings  do define quadrics.)
There is however a
reasonable guess as to the correct answer that goes back to the
translation
principle for simple Lie algebras. For simplicity assume that
$k=\mathbb C$.
Then $\Skl_4$ can be thought of as a deformation of the {\it
homogenized
enveloping algebra} of the complex Lie algebra $\mathfrak{sl}(2)$ (this
is the
ring obtained by taking $\alpha_i=0$ for all $i$ in the definition of
$\Skl_4$). We want to regard $\Skl_4$  as a deformation of
$U(\mathfrak{sl}(2))$ since then the   $R_{\lambda}$ correspond to
the
infinite dimensional  primitive factors $U_\mu=
U(\mathfrak{sl}(2))/(\Omega-\mu)$, where $\Omega$ is the Casimir element
and
$\mu\in k$. The translation principle \cite{BB} asserts that certain
$U_{\mu}$
are Morita equivalent and  by \cite{hodges1} there are no other Morita
equivalences.

  Much of this generalizes. Notably, there exists an
analogous translation principle for $\Skl_4$ (see \cite{VdB11})
 and  for 
noncommutative quadrics (see \cite{VdB26}). The translation principle
for noncommutative quadrics comes out of the observation that, in
contrast to Corollary~\ref{noncanonical},  the sequence
$(\Oscr(n))_n$ is not now determined up to shift. For example, on
$\PP^1\times \PP^1$ we could equally well have chosen
\[
\Oscr(n)=
\begin{cases} 
\Oscr(k,k+a)&\text{if $n=2k$}\\
\Oscr(k,k+1+a)&\text{if $n=2k+1$}
\end{cases}
\]
in place of \eqref{exceplb}.
Analyzing the possible sequences $(\Oscr(n))_n$  leads to
the desired translation principle.
The analogue of
Hodges' theorem, which asserts that there are no other Morita equivalences,
 is also likely to be  true. This has been proved if the triple
 $(Y,\Lscr_0,\Lscr_1,\Lscr_2)$ is sufficiently generic but in the
 general case a few technical problems remain which the second author is
 currently trying to solve.


\section{$\PP^1$-bundles}
\label{sec11}
A question that was not explicitly discussed in the last section is
whether the quadratic quotients $R_{\lambda}$ 
of the Sklyanin algebra $\Skl_4$ do in fact correspond to
noncommutative quadrics. One way to tackle this is through the
concept of a $\PP^1$-bundle.
Recall that commutative quadrics can be written as a $\PP^1$-bundle
$\PP(\Escr)$\label{p1bundle-index} where
$\Escr=\Oscr_{\PP^1}^2$ \cite[\S II.7]{H}.
 In this section  we give similar results for
noncommutative quadrics (see Theorem \ref{another-helix})
 and for the $R_{\lambda}$
 (Proposition \ref{trans}). This shows that, generically,
  the $R_{\lambda}$ are 
 noncommutative quadrics. 

 In order to explain this material we have to introduce some new
concepts. Some of these definitions are a little technical, but the
reader  will
not go far wrong if, initially, he/she thinks of bimodules in the more
mundane, ring theoretic  sense.

Below  $W,X,Y,Z$ are schemes of finite type over $k$. A {\it coherent
$\Oscr_X\hyphenn \Oscr_Y$-bimodule} $\Escr$ is by definition a coherent
$\Oscr_{X\times Y}$-module such that the  support of $\Escr$ is finite
over
both $X$ and $Y$.  We denote the abelian
category of coherent $\Oscr_X\hyphenn \Oscr_Y$-bimodules by
$\bimod(\Oscr_X\hyphenn \Oscr_Y)$\label{OX-bimod-index}.
 Any  coherent bimodule 
$\Escr$ defines a right exact functor $-\otimes_{\Oscr_X} \Escr :
\Qch(X)\r
\Qch(Y)$ via 
 $\pr_{2\ast}(\pr_1^\ast(-)\otimes_{\Oscr_{X\times Y}}\Escr)$. 
Similarly, if $\Fscr$ is a coherent  $\Oscr_Y\hyphenn \Oscr_Z$-bimodule
then the tensor product
$\Escr\otimes_{\Oscr_Y} \Fscr$ is defined as
$\pr_{13\ast}(\pr_{12}^\ast\Escr\otimes_{\Oscr_{X\times Y\times Z}}
\pr_{23}^\ast\Fscr)$. It is easy to show that this definition yields all
the expected properties.

Now  assume that we have finite maps
$\alpha:W\r X$ and  $\beta: W\r Y$. If $\Hscr$ is a coherent
$\Oscr_W$-module then  we  denote the $\Oscr_X\hyphenn \Oscr_Y$-bimodule 
$(\alpha,\beta)_\ast\Hscr$ by ${}_\alpha\Hscr_\beta$.  Any bimodule
$\Escr$ can be presented in this form
since we may take $W$ to be the support of $\Escr$. 
{}From the definition it is easy to check that $-\otimes
{}_\alpha\Hscr_\beta=\beta_\ast(\alpha^\ast(-)\otimes_{\Oscr_W} \Hscr)$.

\begin{example} Let $\Lscr$ be a line bundle on $X$ and let $\sigma\in
  \Aut(X)$. Put $\Mscr={}_\id \Lscr_\sigma$. Then we obtain
  $-\otimes_{\Oscr_X} \Mscr=\sigma_\ast(-\otimes_{\Oscr_X}\Lscr)$.
  Thus $-\otimes_{\Oscr_X}\Mscr$ is the standard autoequivalence of
  $\Qch(X)$ that was introduced in \S\ref{sec2.2} and  played such  
  a fundamental role in the construction of twisted homogeneous
  coordinate rings.
\end{example}

We will define graded
algebras in  $\bimod(\Oscr_X\hyphenn \Oscr_X)$ as  formal direct sums 
$\Ascr=\bigoplus
\Ascr_n$, for $\Ascr_n\in \bimod(\Oscr_X\hyphenn \Oscr_X)$, 
equipped with
multiplication maps $\Ascr_m\otimes_{\Oscr_X} \Ascr_n\r \Ascr_{m+n}$ 
and a unit map
$\Oscr_X\r \Ascr_0$ satisfying the usual properties. A graded right
$\Ascr$ module is a formal direct sum $\Mscr=\bigoplus_n \Mscr_n$, for
$\Mscr_n\in \Qch(X)$,  
equipped with multiplication maps $\Mscr_m\otimes_{\Oscr_X} \Ascr_n\r
\Mscr_{m+n}$,
again satisfying the usual properties. We denote the corresponding
category by $\Gr(\Ascr)$. The quotient category $\QGr(\Ascr)$ is defined
in the expected way, and other constructions like $\ZZ$-algebras
generalize
routinely to this context. The formal definitions can be found in 
\cite{VdB19}.

We now turn to the four dimensional Sklyanin algebra 
$\Skl_4=\Skl_4(E,\Lscr,\sigma)$ and its factors
$R_{\lambda}(E,\Lscr,\sigma)
=\Skl_4/(\Omega_\lambda)$, for $\lambda\in
\PP^1$, as defined in Example~\ref{skl4}. Fix a map $\omega:E\r \PP^1$
of
degree two.

\begin{proposition} \cite[Theorem 7.2.1]{VdB11} \label{trans}
Let $R_{\lambda}=R_{\lambda}(E,\Lscr,\sigma)$   
and assume that $|\sigma|=\infty$.  
Then, for sufficiently
general $\lambda$, we have $\QGr(R_{\lambda})=\QGr(T\Escr/(Q))$ where
$\Escr=
{}_\omega \Vscr _{\omega\sigma}$  for  some line bundle $\Vscr=\Vscr(\lambda)$
 of degree
four on $E$,  
$T\Escr$ is the tensor algebra of $\Escr$ and $Q$  is an invertible
subbimodule of $\Escr\otimes \Escr$.  
\end{proposition}

Since $\Escr$ is locally free of rank two on both sides we may clearly view 
$\QGr(T\Escr/(Q))$  as a  noncommutative analogue of $\PP(\Escr)$. 
We have therefore written $\QGr(R_{\lambda})$ as a noncommutative
analogue of a
$\PP^1$-bundle on $\PP^1$.

An obvious question raised by   
Proposition~\ref{trans} is to determine 
when it is possible to form sheaves of algebras of the form $T\Escr/Q$,
where $\Escr$ is a rank two 
locally free $\Oscr_X$-bimodule $\Escr$ over a smooth curve $X$
and $Q$ is an appropriate subbimodule of $\Escr\otimes \Escr$.
One solution is given  by  Patrick in \cite{Pat1,Pat2}. 
In the commutative case, $Q$ is just the subbimodule of 
$\Escr\otimes \Escr$ generated 
by sections of the form $x\otimes y-y\otimes x$. 
Thus, in general, one at least needs $Q$ to be a rank one subbimodule of 
$\Escr\otimes \Escr$ that cannot be decomposed as 
$Q_1\otimes Q_2$ for $Q_i\subset \Escr$.
This puts 
severe restrictions on the possible choices of~$\Escr$,
to the extent that Patrick can classify them into 
just four basic classes \cite{Pat1}.

This suggests that $T\Escr/Q$ may be the wrong model for $\PP(\Escr)$.
As the second author has recently observed,
 working in the more general setting of $\ZZ$-algebras allows 
one to form an appropriate analogue of $\PP(\Escr)$ without any such
 restrictions on the rank two bimodule $\Escr$.
This even enables one to 
 define the categories $\Gr(T\Escr/(Q))$ and $\QGr(T\Escr/(Q))$ from
 Proposition~\ref{trans} independently of $Q$.
 We briefly sketch the construction.

Assume that $X$ and $Y$ are smooth schemes   and let $\Escr\in
\bimod(\Oscr_X\hyphenn\Oscr_Y)$ be
locally free of rank two on both sides.
Then it is   easy to show that $-\otimes_{\Oscr_X} \Escr$ has
a right adjoint
$-\otimes_{\Oscr_Y}\Escr^\ast$, where 
$\Escr^\ast\in \bimod(\Oscr_Y\hyphenn\Oscr_X)$ is also locally free of
rank two on both sides.
Repeating this construction,  we may define $\Escr^{\ast 2}
=\Escr^{\ast\ast}$ by requiring that 
$-\otimes_{\Oscr_X}\Escr^{\ast\ast} $ is the right adjoint of 
$-\otimes_{\Oscr_Y}\Escr^\ast$. By induction, we define $\Escr^{\ast n}
=(\Escr^{\ast (n-1)})^\ast$ for $n\ge
0$ and by  considering left adjoints we may define
$\Escr^{\ast n}$ for $n<0$. Standard properties of adjoint
functors yield  a bimodule inclusion
$
j_n:\Oscr_{X_n}\hookrightarrow \Escr^{\ast n}\otimes_{\Oscr_{X_{n+1}}}
\Escr^{\ast (n+1)},
$ where $X_{n}$ is respectively $X$ or $Y$ 
depending on whether $n$ is even or odd.

We now define $S\Escr$ as the $\ZZ$-algebra which satisfies
\begin{itemize}
\item[(a)]
 $(S\Escr)_{nn}=\Oscr_{X_n}$;
\item[(b)]
$(S\Escr)_{n,n+1}=\Escr^{\ast n}$;
\item[(c)] 
$S\Escr$ is freely generated by the $(S\Escr)_{n,n+1}$,
  subject to the relations given by the images of $j_n$.
\end{itemize}
It is easy to show:

\begin{proposition}\label{p1bundle3} Let $X=Y=E$ and $\Escr$ be defined by 
Proposition~\ref{trans}.  Then  the  category $\Gr(T\Escr/(Q))$ which
appeared in Proposition~\ref{trans} 
 is equivalent to $\Gr(S\Escr)$ Thus  the
construction is independent of $Q$!
\end{proposition}

Given the success of this approach, we propose the
following definition.

\begin{definition} 
\label{p1bundle} 
{\em A  noncommutative
$\PP^1$ bundle} on a smooth scheme $X$ of finite type
  over $k$ is a category of the form $\QGr(S\Escr)$ where
   $\Escr\in \bimod(\Oscr_X\hyphenn\Oscr_Y)$, for some smooth scheme
$Y$,
   such that $\Escr$ is locally free  of rank two on both sides.
\end{definition}

We leave it to the reader to check the following result.
\begin{proposition} \label{p1bundle2}
Let $\Escr$ be  as in Definition  \ref{p1bundle} and
  assume that $\Uscr,\Vscr$ are line bundles on $X$ and $Y$
  respectively. Put
  $\Fscr=\Uscr\otimes_{\Oscr_X}\Escr\otimes_{\Oscr_Y} \Vscr$. Then
  $\Gr(S\Escr)$ and $\Gr(S\Fscr)$ are equivalent.
\end{proposition}

Finally we mention the following result from \cite{VdB26}. 

\begin{theorem} \label{another-helix}
Assume that $A$ is a $\ZZ$-algebra satisfying
  (Q1)-(Q4) from \eqref{noncom-quadrics} 
  with associated quadruple $(Y,\Lscr_0,\Lscr_1,\Lscr_2)$.
Assume that, in the corresponding sequence  $(Y,(\Lscr_i)_i)$ defined by 
 \eqref{sequence}, we
have $\Lscr_0\not\cong \Lscr_{2n+1}$ for $n\ge 0$. Let $Y'$ be the
image of $Y$ in $\PP^1\times \PP^1$ under the map defined by
$(\Lscr_0,\Lscr_2)$ and let $\Lscr_{1}'$ be the direct image of
$\Lscr_{1}$ under the induced map $Y\r Y'$.  
Then $\QGr(A)=\QGr(S\Mscr)$ for $\Mscr={}_{\pr_1}(\Lscr_1')_{\pr_2}$.  
\end{theorem}
Combining (\ref{trans}--\ref{another-helix})
shows that {\it  the  quotients $R_\lambda$
of the four-dimensional Sklyanin algebras are noncommutative quadrics.}

\section{Noncommutative Blowing Up} \label{sec12}

In this section we summarize some recent results on blowing up of
noncommutative surfaces. In particular, this allows us to give
noncommutative
analogues of the construction  of cubic surfaces, Cremona transformation
and
quadrics that were outlined in \S\S\ref{seccub}--\ref{secquad}, as well
as
answer the question posed by Example~\ref{moredetailsfollow}. This
material is
all taken from \cite{VdB19, VdB27}, to which the reader is referred for
more
details.

\subsection{Motivation, definition and properties}
\label{sec12.1}
To fix ideas let us first consider the affine case. Assume that $R$ is a 
finitely generated  noetherian  $k$-algebra and let $\mathfrak m$ be a  
maximal ideal of $R$ such that $R/\mathfrak m=k$. 
Then, as in \eqref{standardreesring},
 a natural idea is to
define the blow-up of $R$ at $\mathfrak m$ as  $\QGr(D)$  where 
$D=R\oplus \mathfrak m\oplus \mathfrak m^2\oplus \mathfrak
m^3\oplus\cdots$ is
the Rees  algebra  associated to $\mathfrak m$.
Unfortunately this na{\"\i}ve definition is often faulty, as is shown by 
the following example
from \cite{Ar2}. 

Let $R=k\{ x,y\}/(yx-xy-y)$ with 
$\mathfrak m=(x,y)$.  Thus  $R$ is the enveloping algebra of the
two-dimensional 
solvable Lie algebra and resembles the polynomial ring in two variables
that
 was considered in \S\ref{seccomsurov}.
However,  $\mathfrak m^n=(x^n,y)$. Hence the analogue of the
exceptional curve
\[
D/\mathfrak mD=R/\mathfrak m\oplus 
\mathfrak m/\mathfrak m^2\oplus \mathfrak m^2/\mathfrak m^3\oplus\cdots
\]
is isomorphic to $k[x]$. Thus $\QGr(D/\mathfrak m)$ is a point, whereas
intuitively we would expect it to be one-dimensional.

What has gone wrong is that, unlike the commutative case, $\mathfrak m$
is ``linked'' to other prime ideals in the sense that there are 
other simple modules $R/\mathfrak m_i$ that appear in the injective hull 
$E(_RR/\mathfrak m)$. Somehow, we need to involve them in the definition
of
blowing-up. A general discussion of linked primes can be found in
\cite[Chapters~11,12]{GW} (see, in particular \cite[Example 11G]{GW} for
a
discussion of the ring $R$), but we will use a simpler approach for this
example. The basic observation is that the commutator ideal of $R$ is 
$yR$ and the element $y$  induces an automorphism $\tau$ of $R$ 
by $a\mapsto y^{-1}a y$. 
The prime ideals $R/\mathfrak m_i$ mentioned above are exactly the
$\tau^i(\mathfrak m)$. The correct  blow-up of $R$ at $\mathfrak m$
turns out
to be $\QGr(D)$, for  
\[
D=R\oplus I_1\oplus I_2\oplus I_3\oplus\cdots\qquad \text{where}\quad
I_n =\mathfrak m\tau(\mathfrak m)\cdots
\tau^{n-1}(\mathfrak m)_{\tau^n}
\]
(thus, the right action of $R$ on $I_n$ is twisted by $\tau^n$).
 If $L=D/\tau^{-1}(\mathfrak m) D$, then
 $\dim L_u=u+1$ and so $L$ plays the role of the exceptional curve.
We refer the reader to  \cite{Ar2} for a detailed workout of
this example.

This construction generalizes to any 
 ring whose commutator ideal is invertible.
The latter hypotheses is not unreasonable since, if we look at the case
of
a Poisson surface,   we see that we expect a noncommutative smooth
surface to contain a commutative curve.

However, we also want   blow-ups to apply in more general
situations where the concept of the commutator ideal does not really
make sense, in particular since we want to blow up at one-sided ideals. 
One such situation is given by the noncommutative
$\PP^2$'s and $\PP^1\times\PP^1$'s which were introduced in
\S\ref{secncpq}. For simplicity let us consider the nontrivial 
noncommutative $\PP^2$'s, as introduced in Definition~\ref{noncom-p2}.
By construction, such a category $\QGr(A)$ contains a subcategory 
isomorphic to $\Qch(E)$, for some elliptic curve $E$
 and so, as in \S \ref{secasregular}, we consider $E$ to be a closed
(commutative) subscheme of the noncommutative scheme $\QGr(A)$.

Now let $p\in E$.   The previous discussion suggests that it should be
possible to blow up $p$. However it is not clear how to proceed. As
explained in \S\ref{secasregular}, $p$ corresponds to a  
point module  over $A$ which, since it  only a
right module, cannot be used to construct a Rees algebra.

Such difficulties can be solved by working with categories rather than
algebras.
To explain this, we need a little more   geometric
language.
As before, we define a noncommutative scheme to be a
Grothendieck category but to stress the geometric viewpoint we will
now use roman capitals $X,Y,\ldots$.  If $X$ is a commutative
scheme then we will identify $X$ with $\Qch(X)$. This is justified by
the reconstruction theorem proved by Rosenberg \cite{rosenberg1}.
A morphism $\alpha:X\r Y$ between noncommutative schemes will be a right
exact
functor $\alpha^\ast:Y\r X$ possessing a right adjoint
(denoted by $\alpha_\ast$). In this way the noncommutative schemes form
a
category (more precisely a two-category, see \cite{KS}).

Let  $X,Y,Z$, \dots be noncommutative schemes.  Inspired by
\S\ref{sec11}
we would like to introduce $X\hyphenn Y$-bimodules. Unfortunately the
definition given in \S\ref{sec11} does not generalize to the
noncommutative 
case.
Therefore we change the definition of bimodules: An \emph{$X$\hyphenn$Y$
  bimodule} is just an object in 
$$
  \Bimod(X\hyphenn Y)=\{\text{morphisms $X\r Y$}\}
$$
 We write the functor corresponding to $\Mscr\in
  \Bimod(X\hyphenn Y)$ as $-\otimes_X \Mscr$. 
Likewise we write composition of functors
\[
\Bimod(X\hyphenn Y)\times \Bimod(Y\hyphenn Z)\r
\Bimod(X\hyphenn Z)
\]
as $-\otimes_Y-$. 
It is clear that $\Bimod(X\hyphenn X)$ contains a canonical object
 $o_X$ given by the identity functor. 

We will  need to talk about images and subobjects in
$\Bimod(X\hyphenn Y)$, which is a problem since 
$\Bimod(X\hyphenn Y)$ is  not an
abelian category  (at least not in any obvious way).  We  resolve the
problem
by  embedding  $\Bimod(X\hyphenn Y)$ into the opposite category to the
category of left exact functors from $Y$ to $X$.
As this bigger category is abelian, images and subobjects do exist
there.

As $\Bimod(X\hyphenn X)$ is a monoidal category we can define graded
algebras and graded modules over it in formally the same way as this was
done
in \S\ref{sec11}.

The definition of blowing-up given in \cite{VdB19} works in the
situation where
$X$ contains a commutative curve as a \emph{divisor}. To make this more
precise, we consider  the identity functor $o_X$ on $X$, as defined
above.
 It is clear that $o_X$ is  an algebra and, furthermore, that
$\Mod(o_X)=X$. We will assume that $o_X$ contains an invertible
subbimodule
$o_X(\hyphenn Y)$ (an {\it invertible ideal}) such that
$\Mod(o_X/o_X(\hyphenn
Y))$ is equivalent to $\Qch(Y)$, for some commutative curve $Y$.
We also need some sort of smoothness condition on $X$. Since it is
obviously sufficient to impose this in a neighbourhood of $Y$, we assume
that every object in $\Qch(Y)$ has finite injective dimension in
$X$. Finally we will assume that $X$ is \emph{noetherian}, by which we
mean that $X$ is a locally noetherian category.

We emphasize that all these assumptions hold when $X=\QGr(A)$ is a
nontrivial  noncommutative $\PP^2$, with $Y=E$ the corresponding
elliptic
curve.

Let $p\in Y$. Then $p$ defines a subbimodule $m_p$ of $o_X$ which is
the analogue of the maximal ideal corresponding to $p$ in $\Oscr_X$, if
$X$ were commutative.  We put
$I=m_po_X(Y)$ and define
\[
 \Dscr=o_X\oplus I\oplus I^2\oplus \cdots
 \]
where $I^{n}$ is defined as the image of $I^{\otimes n}$ in
 $o_X(nY)\overset{\text{def}}{=} o_X(-Y)^{\otimes -n}$.  
 This is the Rees algebra associated to $I$ and we define the
 blow-up $\widetilde{X}$ of $X$ in $p$ to be $\QGr(\Dscr)$.

A large part of \cite{VdB19} is devoted to proving that $\widetilde{X}$
satisfies   properties  similar to those of $X$ and, furthermore, that
we have
obtained an analogue of the blow-up of a commutative surface.
In particular, that paper constructs a
commutative diagram of noncommutative noetherian schemes.
\[
\begin{CD}
\widetilde{Y} @>i>> \widetilde{X}\\
@V\beta VV @VV\alpha V\\
Y @>>j> X
\end{CD}
\]
where the horizontal arrows are inclusions, $\widetilde{X}$ is again a
noncommutative noetherian scheme and  $\widetilde{Y}$
is a commutative curve which plays the role of the strict transform of
$Y$. Moreover, $\widetilde{Y}$ is again a divisor in 
$\widetilde{X}$  and every object in
$\Qch(\widetilde{Y})$ has finite injective dimension in $\widetilde{X}$.

\subsection{More general birational transformations}
As indicated in \S\ref{secsomref} it is convenient 
to have available more general
birational transformations that accomplish combinations of
blowing-ups and blowing-downs. It turns out that we can mimic the
constructions given in \S\ref{secsomref} in the  generality of
\S\ref{sec12.1}.

First assume that $d$ is a divisor in $o_Y$ and let $m_{Y,d}$ be the
corresponding ideal in $o_Y$. Let $m_d$ be the inverse image of
$m_{Y,d}$ in $o_X$. If we assume that   we have a sequence 
$(\Lscr_n)_n$ of noetherian objects on $X$, 
then we can define $\widetilde{X}$ as $\QGr(D)$,
where $D$ is the $\ZZ$-algebra given by
\[
D_{a,b}=\Hom_X(\Lscr_{-b},\Lscr_{-a}m_{\tau^{-a}d}\cdots
m_{\tau^{-b+1}d})
\]
In this generality, presumably not much can be said about the
properties of $\widetilde{X}$. However we can use these definitions
 to generalize
the three classical examples given in
\S\ref{seccub}--\ref{secquad}. 

\subsection{The non-commu\-tative Cremona
transformation}\label{noncom-crem}
Suppose that $X$ is a nontrivial  
noncommutative $\PP^2$ (see Definition~\ref{noncom-p2})
and let $Y=E\subset X$ be the
associated  curve. Since $X=\QGr(A)$
for a graded ring $A$, there exist naturally defined objects $\Oscr_X(n)
\in X$, coming from the shifts $A(n)$.
 We put $\Lscr_n=\Oscr_X(2n)$.

Choose three distinct non-collinear points $p,q,r\in Y$ and set
 $d=p+q+r$. Then it turns out that $\widetilde{X}$ as defined
above is also a noncommutative $\PP^2$. Using this we can show that
there are exactly four different types of division algebras that can
occur as the function field of a noncommutative $\PP^2$, namely:  
\begin{itemize}
\item[(a)] Division rings finite dimensional over a commutative field of
transcendence degree two.
\item[(b)] The quotient field of the first Weyl algebra.
\item[(c)] The function fields of the affine quantum planes (see Example
  \ref{moredetailsfollow}). 
\item[(d)] The function fields of the three dimensional Sklyanin
  algebras.
\end{itemize}
Since both 
(b) and (c) are special cases division rings of skew polynomial rings,
this confirms  the conjecture of Artin from \S\ref{birat-class}
in this case.


\subsection{Noncommutative quadrics}\label{noncom-quad}
\label{secnoncomquad}
Now let $X$ be a noncommutative quadric as defined in
\S\ref{sec10}. There exist naturally defined objects $\Oscr_X(n,n)$
(analogues of $\Oscr_{\PP^1}(n)\boxtimes \Oscr_{\PP^1}(n)$)
and we put $\Lscr_n=\Oscr_X(n,n)$.

One can show that the elliptic curve $Y$ used in the construction of
$X$ actually embeds as a divisor in $X$. Let $d=p$ be a point
in $Y$. Then $\widetilde{X}$ turns out to be a 
noncommutative $\PP^2$. From
this one deduces the following result.
\begin{theorems} Every noncommutative quadric is birational to a
  noncommut\-ative $\PP^2$.
\end{theorems}
One can prove the converse of this theorem using similar methods.

\subsection{Noncommutative cubic surfaces}\label{noncom-cub}
We again let $X$  be a nontrivial noncommutative $\PP^2$ but now we put  
$\Lscr_n=\Oscr_X(3n)$ and we take six points $\{p_1,\ldots,p_6\}\in Y$
in
sufficiently general position. Put $d=\sum_i p_i$. Then it follows
from  \cite{VdB19}  that  there exists
 a graded ring $S$, with
the Hilbert series of a polynomial ring in 4 variables, and containing a
normal non-zero divisor 
 $g$ of degree three such that $\widetilde{X}=\QGr
(S/(g))$. In other words,
 $\widetilde{X}$ can be considered as a noncommutative cubic surface
in a noncommutative $\PP^3$. 

\vskip 5pt
In each of these examples (\ref{noncom-crem})--(\ref{noncom-cub}),
the analogy with the commutative
case is striking. 

\subsection{More details on Example \ref{moredetailsfollow}}
\label{secmoredetails} We end the paper by showing how to use the
techniques
of this section to prove the result asserted  in that example.
Recall that 
 $E=k(u,v;vu=quv)$ and $D=E^\tau$ where $\tau u=u^{-1}$, $\tau
v=v^{-1}$. The problem is to determine where $D$ fits in Artin's list of
division rings from \S\ref{birat-class}.

Put $L=k(u)$  and let $K$ be the subfield of $L$ given by
 $k(s)=k(u+u^{-1})$.
Define $\psi\in \Aut(L)$ by
$\psi(u)=\theta u$ where $\theta$ is a scalar satisfying
 $\theta^{-2}=q$. Let $M$ be the $K$-$K$ bimodule given by ${}_\id
 L_\psi$. Thus $s$ acts on the right of $M$ by $t=\theta u+\theta^{-1}
 u^{-1}$.  
An explicit computation now shows that $D$ is the function field of the
$\PP^1$-bundle over $\Spec K$ associated to $M$.

We may identify
$M$ with the $k(s)\otimes_k k(t)$-module $k(u)$ where $s$ and $t$ act
in the  way indicated above. We will now view $K$ as the function
field of $\PP^1$ and we will identify $M$ with the localization at the
generic point of a sheaf-bimodule $\Mscr$. To do this note that the
images of
$s$ and $t$ in $k(u)$ satisfy the equation
\begin{equation}
\label{curveequation}
(\theta s-t)(\theta^{-1} s-t)=-(\theta-\theta^{-1})^2
\end{equation}
Let $Y\subset \PP^1\times \PP^1$ be the closure of the curve in
$\Spec k[s,t]$ defined by \eqref{curveequation}. Closer inspection
reveals that $Y$ is in fact a nodal elliptic curve with a unique
singular point at $(\infty,\infty)$.
 
We now define $\Mscr={}_{\pr_1}\Lscr{}_{\pr_2}$ where $\Lscr$ is a
generic line bundle of degree two on $Y$. Thus $D$ is the function
field of the $\PP^1$-bundle over $\PP^1$ associated to $\Mscr$. By
Theorem~\ref{another-helix} this $\PP^1$-bundle is isomorphic to a
noncommutative quadric and, by \S\ref{secnoncomquad} above, a
noncommutative
quadric is birational to a noncommutative plane $X$. If we examine more
closely the birational equivalence in \S\ref{secnoncomquad} then 
this noncommutative plane happens to be associated to a triple
$(Y,\Lscr,\sigma)$ where $Y$ is a triangle.  Thus the function field
of $X$ is the same as that of an affine quantum plane and we are done.

Of course, once one knows that $D$ is the function field of a quantum plane, 
it should not be too difficult to write down the requisite generators
and recently we have succeeded in achieving this. 
Set 
$$f=(uv-u^{-1}v^{-1})(v-v^{-1})^{-1}\qquad{\rm and }\qquad 
g=(u-u^{-1})(v-v^{-1})^{-1}.$$
Then one can show that $fg=qgf$ and that $D$ is the division ring of fractions
of the quantum plane $k\langle f,g\rangle$. In particular, $D\cong E$.
 Proving this is a long but 
elementary computation.

\clearpage
\section*{Index}\label{Index}

\noindent
Much of the basic notation can be found in the Notation Section~\ref{notation}
and concepts defined there are not repeated in this index.

\begin{multicols}{2}

ample divisor, \pageref{amplediv-index} 
 
ample pair $(\Oscr,s)$, \pageref{ampleness}

AS Gorenstein, \pageref{asg-index} 

AS regular, \pageref{AS-regular}

$\bimod(\Oscr_X\hyphenn \Oscr_Y)$, \pageref{OX-bimod-index}

birational equivalence, \pageref{birat-index}

$\chi$ condition, \pageref{chi}

  classical Serre duality, \pageref{serredual-index}

cohomological dimension, \pageref{cohomdim}

exceptional curve, \pageref{excep-curve}

exceptional object, \pageref{excep-ob-index}

 $\Ext$-finite, \pageref{extfin-index}
 
finite length object, \pageref{finlengthob-index}

 Gelfand-Kirillov dimension, \pageref{GKdim-index}
 
Goldie ring, \pageref{goldie-index} 
 
  hereditary category, \pageref{heredcat-index}
  
   homological dimension $n$, \pageref{homdim-index}

  maximal order, \pageref{maxorder-index}

  nice ring, \pageref{nice-index}
 
  noncommutative $\PP^1$-bundle, \pageref{p1bundle}

noncommutative quadric, \pageref{noncom-quadrics}

(nontrivial) noncommutative $\PP^2$, \pageref{noncom-p22}

    one-critical object, \pageref{1critob-index} 

$\PP^1$-bundle $\PP(\Escr)$, \pageref{p1bundle-index}

$\pi(M)$, \pageref{pi-index}

point scheme  $\Pscr(A)$, \pageref{pointscheme-index}

 quasi-unipotent automorphism, \pageref{quasiunipot-index}
 
$R_{\lambda}$, \pageref{Rlambda-index}

saturated, \pageref{saturated-index}

Serre duality, Serre functor, \pageref{serre-fn}

$\sigma$-ample divisor, \pageref{s-ample-index}

Sklyanin algebras $\Skl_3$, $\Skl_4$, \pageref{skl}, \pageref{skl4}

  twisted  coordinate   ring, \pageref{sec2.2}
  
  twisting ring, \pageref{twistingring-index} 

Veronese ring, \pageref{veronese-index}

$\ZZ$-algebra, \pageref{Z-algebra-index}

\end{multicols}

\ifx\undefined\bysame
\newcommand{\bysame}{\leavevmode\hbox to3em{\hrulefill}\,}
\fi

 \end{document}